\documentclass[12pt,a4paper,twoside, english]{amsart}
\usepackage[a4paper, left=25mm, right=25mm, top=25mm, bottom=26mm]{geometry}
\usepackage{amsmath,amsfonts,amssymb}
\usepackage{float}
\usepackage{babel}
\usepackage{hyperref}

\usepackage{wrapfig}
\usepackage{booktabs}

\usepackage{stmaryrd}
\usepackage{comment}
\usepackage{caption}
\usepackage{sidecap}
\usepackage{xcolor}

\usepackage{graphicx}
\usepackage{subcaption}

\hypersetup {colorlinks=true, linkcolor=blue, urlcolor=blue, citecolor=blue}

\usepackage{lmodern}
\usepackage{microtype}

\usepackage[
  backend=biber,
  style=numeric,        
  sorting=nyt,
  maxbibnames=5,        
  minbibnames=1,
  giveninits=true,      
  doi=true,
  url=false,
  isbn=false,
  eprint=false
]{biblatex}

\newcommand\R{\mathbb{R}}
\newcommand\C{\mathbb{C}}
\newcommand\Z{\mathbb{Z}}
\newcommand\eps{\varepsilon}
\DeclareMathOperator{\curl}{\mathrm{curl}}

\newcommand\cC{\mathcal{C}}

\newcommand\cH{\mathcal{H}}

\newcommand\cZ{\mathcal{Z}}

\newcommand{\sign}{\operatorname{sign}}

\newcommand{\del}{\partial}
\newcommand{\ra}{\rangle}
\newcommand{\la}{\langle}

\newcommand{\one}[1]{{\bf{1}}_{#1}}

\newcommand{\loc}{\mathrm{loc}}

\newcommand\vth{\vartheta}
\newcommand\fhi{\varphi}
\newcommand{\weakto}{\rightharpoonup}

\usepackage{mathtools}
\usepackage{extarrows}
\renewcommand{\div}{\operatorname{div}}

\newtheorem{theorem}{Theorem}[section]
\newtheorem{definition}[theorem]{Definition} 
\newtheorem{lemma}[theorem]{Lemma}

\newtheorem{remark}[theorem]{Remark}

\def\XXint#1#2#3{{\setbox0=\hbox{$#1{#2#3}{\int}$ }
		\vcenter{\hbox{$#2#3$ }}\kern-.6\wd0}}

\numberwithin{equation}{section}

\allowdisplaybreaks

\newcommand{\Oe}{{\Omega_\eta}}

\newcommand{\Ee}{E^\eta}
\newcommand{\He}{H^\eta}
\newcommand{\rot}{\operatorname{curl}}

\begin{document}
	
\pagestyle{myheadings} \markboth{Maxwell's equations with thin wires
  in critical scaling}{B.~Schweizer and D.~Wiedemann}
	
\thispagestyle{empty}
\begin{center}
  ~\vskip3mm {\Large\bf Maxwell's equations with thin \\[2mm]
    wires in critical scaling}\\[5mm]
  {\large B.~Schweizer\footnotemark[1] and
    D.~Wiedemann\footnotemark[1]}\\[3mm]
  
  13.~November 2026 \\[3mm]

\end{center}

\footnotetext[1]{Technische Universität Dortmund, Fakult\"at f\"ur
  Mathematik, Vogelpothsweg 87, D-44227 Dortmund,
  ben.schweizer$@$tu-dortmund.de, david.wiedemann$@$tu-dortmund.de}
	
\begin{center}
  \vskip4mm
  \begin{minipage}[c]{0.87\textwidth}
    {\bf Abstract:} We derive effective equations for a situation in
    which electromagnetic waves hit an array of very thin wires. We
    study two cases, perfectly conducting wires and highly conducting
    wires, and investigate in both the critical scaling of the
    parameters. The electric field $\Ee$ and the magnetic field $\He$
    satisfy the time-harmonic Maxwell equations, and the small
    parameter $\eta>0$ is the spacing of the wires. A non-critical
    scaling of the wires implies either that the wires have no effect
    or that the tangential component of the limiting electric field in
    the direction of the wires vanishes along the effective
    interface. In the critical scaling, we observe an interesting
    effect: The wires act as oscillators and there is a relation
    between the electric field in longitudinal direction and the jump
    of the magnetic field in transversal direction involving second
    derivatives of a flux quantity.


    

  \end{minipage}\\[6mm]
\end{center}

\nocite{Delourme-Schweizer-Wiedemann-Pol-2025}

\section{Introduction}

Maxwell's equations describe electromagnetic waves with an electric
field $E$ and a magnetic field $H$.  Since we are interested in
complex geometries, we consider domains
$\Oe \subset \Omega\subset\R^3$ that depends on a (small) parameter
$\eta>0$. We are interested in the situation that many thin wires
interact with the electromagnetic waves; we will always assume that
the wires occupy a domain $\Sigma_\eta\subset\Omega$ and consider the
domain $\Omega_\eta \coloneqq \Omega\setminus \overline{\Sigma_\eta}$.

As the domain depends on $\eta$, also the solutions depend on $\eta$;
we write $\Ee$ and $\He$.  In the time-harmonic setting with a fixed
frequency $\omega>0$ and with fixed permeability $\mu$ and fixed
permittivity $\eps_0$, the system reads
\begin{subequations}\label{eq:Maxwell:eta:pc}
  \begin{align}\label{eq:Maxwell:eta:pc:1}
    \rot\Ee &= i \omega \mu \He + f_h &&\textrm{ in } \Oe \, , 
    \\\label{eq:Maxwell:eta:pc:2}
    \rot\He &= - i \omega \eps_0 \Ee +f_e &&\textrm{ in } \Oe \, ,
    \\\label{eq:Maxwell:eta:pc:3}
    \Ee \times \nu &= 0 &&\textrm{ on } \partial \Oe \, .
  \end{align}
\end{subequations}
The last equation is a boundary condition that models a perfectly
conducting material outside $\Omega_\eta$; it uses the exterior normal
$\nu$ of $\Omega_\eta$. It demands, in particular, that the wires are
perfectly conducting.  The index $\eta > 0$ indicates that we analyze
a fine-scale structure. The solution is a map
$(\Ee,\He) \colon \Oe\ni x\mapsto (\Ee,\He)(x) \in \C^3\times \C^3$.
The right-hand side $f_e = f_e(x)$ of the second equation models
prescribed external currents; we include $f_h = f_h(x)$ to have the
system more symmetric. In order to simplify calculations, we always
assume that the support of the sources has a positive distance from
the interface $\Gamma$.

\smallskip In the second case of highly conducting materials in the wires,
we consider solutions
$(\Ee,\He) \colon \Omega\ni x\mapsto (\Ee,\He)(x) \in \C^3\times \C^3$
and the system involves a permittivity $\eps_\eta$ that depends on
$\eta$:
\begin{subequations}\label{eq:Strong:Maxwell:eta}
  \begin{align}\label{eq:Strong:Maxwell:eta:1}
    \rot\Ee &= i \omega \mu \He + f_h &&\textrm{ in } \Omega \, , 
    \\\label{eq:Strong:Maxwell:eta:2}
    \rot\He &= - i \omega \eps_\eta \Ee +f_e &&\textrm{ in } \Omega \, ,
    \\\label{eq:Strong:Maxwell:eta:3}
    \Ee \times \nu &= 0 &&\textrm{ on } \partial \Omega \, .
  \end{align}
\end{subequations}
The system with highly conducting thin wires becomes interesting when
the conductivity is large. We set
\begin{equation}
  \label{eq:eps-eta}
  \eps_\eta(x) \coloneqq
  \begin{cases}
    \eps_0\qquad &\text{ for } x\in \Omega_\eta\,,\\
    \eps_0 + \frac{i\kappa}{\theta_\eta}\qquad &\text{ for } x\in \Sigma_\eta\,,
  \end{cases}
\end{equation}
where $\theta_\eta>0$ is a parameter with $\theta_\eta\to 0$ for
$\eta\to 0$. We use $\kappa>0$ as a parameter that characterizes the
conductivity of the wire.

\subsection{The effective system}

We are interested in a system of equations that describes
electromagnetic waves for small $\eta>0$. The precise task is the
following: Assume that $(E,H)$ is a limiting solution in the sense
that, along a subsequence,
\begin{equation}
  \label{eq:weak-limits}
  E^\eta\weakto E\quad\text{ and }\quad H^\eta\weakto H\quad
  \text{ in } L^2(\Omega,\C^3)\,.
\end{equation}
We want to know the equations that characterize the limit $(E,H)$.
When the wires concentrate at an interface $\Gamma$ and when the wires
are all pointing in direction $e_1$, then, in the strong form, we
derive the following limit problem:
\begin{subequations}\label{eq:Strong:Maxwell:limit}
  \begin{align}\label{eq:Strong:Maxwell:limit:1}
    \rot E &= i \omega \mu H + f_h
    &&\textrm{ in } \Omega \, , 
    \\\label{eq:Strong:Maxwell:limit:2}
    \rot H &= - i \omega \eps_0 E +f_e
    &&\textrm{ in } \Omega \setminus \Gamma \, ,
    \\
    \label{eq:Strong:Maxwell:limit:3}
    E \times \nu &= 0
    &&\textrm{ on } \partial \Omega \,,\\
    \label{eq:Strong:Maxwell:limit:int-triv}
    \llbracket H_1 \rrbracket_\Gamma &= 0
    &&\textrm{ on } \Gamma\,,\\
    \label{eq:Strong:Maxwell:limit:int}
    &\text{an interface relation coupling $E_1$ and 
    $\llbracket H_2 \rrbracket_\Gamma$ }
           &&\textrm{ on } \Gamma\,,
  \end{align}
\end{subequations}
where jumps such as $\llbracket H_2 \rrbracket_\Gamma$ are defined as
the trace (in a distributional sense) on $\Gamma$ taken from above
(from $\Omega_+$) minus the trace on $\Gamma$ taken from below (from
$\Omega_-$).  Obviously, the interesting part is the interface
condition \eqref {eq:Strong:Maxwell:limit:int}.  Equations \eqref
{eq:Strong:Maxwell:limit:1}--\eqref{eq:Strong:Maxwell:limit:3} are
immediate consequences of the equations for $(E^\eta,H^\eta)$ and the
weak convergence \eqref {eq:weak-limits}. The interface condition
\eqref{eq:Strong:Maxwell:limit:int-triv} was obtained in
\cite{Schweizer-Wiedemann-Pol-2025}: Since the wires are disconnected
in direction $e_2$, the field $H_1$ cannot jump across $\Gamma$.

In our main result, we determine the transmission condition
\eqref{eq:Strong:Maxwell:limit:int} in the critical scaling. 

\smallskip
\paragraph{\em Critical scaling.} A non-trivial coupling appears in
\eqref{eq:Strong:Maxwell:limit:int} when the wires have a critical
scaling. We consider wires of the relative size $\delta_\eta$, using a
positive parameter $\alpha$ in the definition of $\delta_\eta$.  For
highly conducting wires, the conductivity is scaled by a parameter
$\theta_\eta$:
\begin{equation}
  \delta_\eta =  e^{-\alpha\pi/\eta}\qquad\text{ and }\qquad
  \theta_\eta = \eta\delta_\eta^2\,. 
  \label{eq:delta-theta}
\end{equation}
We denote the two-dimensional volume of the cross-section $\cC_{2D}$
by $\beta \coloneqq |\cC_{2D}|$.  As a consequence, the scaled
cross-section satisfies $|\cC_\eta| = \beta \delta_\eta^2$.

\smallskip
\paragraph{\em Discussion of the interface condition.}

In the case that the scaling is not critical, the limit law was
determined in \cite{Schweizer-Wiedemann-Pol-2025} and
\cite{DelourmeHewett2020} for perfect conductors. The first case is
that the wires are too thin. In this case, they are invisible in the
limit and the effective law is simply
$\llbracket H_2 \rrbracket_\Gamma = 0$.

The second case is that the wires are sufficiently thick. In this
case, the effective law is $E_1|_\Gamma = 0$. This coincides with our
intuition: As the wires are perfectly conducting and $e_1$ is a
tangential vector, there holds $E_1^\eta|_{\del\Sigma_\eta} = 0$. In
some sense, sufficiently thick wires imply that the effective field
also vanishes along the limiting surface $\Gamma$. We note that
\cite{Schweizer-Wiedemann-Pol-2025} allows also wires that are
non-uniform in $e_1$-direction such as bent wires or wires with slits.

Based on the results described above and related problems on the
Poisson \cite{MR1493040} and Stokes problem \cite{MR1079190}, one
might expect that, in the critical scaling of the wires, the result is
a limit law of the form
$\llbracket H_2 \rrbracket_\Gamma = \lambda\, E_1|_\Gamma$ for some
parameter $\lambda\in \C$. This is not the case! Instead, a second
derivative in the wire direction appears. We obtain for \eqref
{eq:Strong:Maxwell:limit:int} in the case of high conductors
\begin{equation}
  \label {eq:int-law-hc-intro}
  \left(\omega^2 \mu\eps_0 + \del_1^2
    + \frac{2i \eps_0}{\alpha\beta  \kappa} \right) j
  =  -\frac{2\eps_0}\alpha  E_1|_\Gamma
  \quad\text{ with }\quad
  j = \frac{1}{i \omega} \llbracket H_2 \rrbracket_\Gamma\,,
\end{equation}
compare \eqref {eq:int-law-hc}.  We note that this limit law admits
resonances in the wires due to the longitudinal second derivative. We
mention that the additional term cannot appear in a strictly
two-dimensional analysis; accordingly, it does not appear, e.g., in
\cite {hewett2016faraday}.

A similar term with second derivatives appears in the seminal paper
\cite{MR2262964} of G.~Bouchitt\'e and D.~Felbacq. They consider thin
parallel wires with circular cross section that fill the whole volume.
The astonishing effective law (3.1) in \cite{MR2262964} reads (we take
the liberty to replace $3$ by $1$ as their wires are pointing in
direction $e_3$ and ours in direction $e_1$ and replace $\kappa$ by
$\kappa_{\text{BF}}$ for clearer distinction):
\begin{equation}\label{eq:Helmholtz-BF}
  \del_1^2 J + \left(k_0^2 + \frac{2 i \pi \gamma}{\kappa_{\text{BF}}} \right) J
  = 2 i \pi \gamma E_1\,.
\end{equation}
We emphasize that \eqref {eq:Helmholtz-BF} describes a law in the
bulk.

The flux quantity $J$ of \cite{MR2262964} is related to our flux $j$
with a complex factor, $j = i \eps_0 J$. Using this identity and
$k_0^2 = \omega^2 \mu\eps_0$, we see that the two formulas \eqref
{eq:int-law-hc-intro} and \eqref {eq:Helmholtz-BF} coincide. Let us
perform the comparison in detail: Their parameter $\gamma$ models the
exponential scaling of the radius of the wires, thus
$\gamma^{-1} = \alpha\pi$. This already shows that the right-hand
sides of \eqref {eq:int-law-hc-intro} and \eqref {eq:Helmholtz-BF}
(the latter with the factor $i \eps_0$ to match $J$ with $j$)
coincide.

For circular wires, there holds $\beta = \pi$. It therefore remains to
compare our $\alpha\beta\kappa/\eps_0 = \alpha\pi\kappa/\eps_0$ with
their $\kappa_{\text{BF}}/(\pi\gamma) = \alpha\,\kappa_{\text{BF}}$.
We can indeed identify $\kappa_{\text{BF}} = \kappa \pi/\eps_0$ -- for
the factor $\eps_0^{-1}$ see \cite[equation (1.2)]{MR2262964} and, for
the factor $\beta$, see \cite[equations (2.2)--(2.3)]{MR2262964}, where
the factor $\pi$ is included in $\kappa_{\text{BF}}$ via
$\theta_\eta$.

We emphasize that the scaling of $\theta_\eta$ is different in the two
works as the volume filling case requires smaller obstacles. This is
just as in more classical results, see \cite{MR1493040} for the
Poisson equation and \cite{MR1079189, MR1079190} for the Stokes
equations.

Loosely speaking, we may say: We obtain the law of \cite{MR2262964}
along the interface, even though the law was derived originally for
wires that fill the bulk.  Perfect conductors can be described by
sending $\kappa\to \infty$; indeed, for perfect conductors, the third
term in the bracket in \eqref {eq:int-law-hc-intro} is not present.

\subsection{Further literature}
A classical analysis of critically scaled inclusions with the famous
``strange term coming from nowhere'' was given in \cite{MR1493040} for
the Poisson problem. This creation of additional terms in the
effective equations was observed for other systems as well: The
Navier--Stokes equations with critically scaled inclusions lead to a
similar zero order term, which is called Brinkman term
\cite{MR1079189, MR1079190}. These works address obstacles that are
either distributed in the bulk or at an interface.  Critically scaled
inclusions distributed in the bulk for Maxwell's equations are
addressed in \cite{MR2262964}. The $\curl$--$\curl$ structure of the
Maxwell equations leads to additional difficulties in the
homogenization.  Resonances in effective Maxwell-system can be used in
multiple ways; we refer to negative index materials in \cite
{Bouchitte-Schweizer-Max-2010} and to perfect transmission effects in
\cite {Bouchitte-Schweizer-Plasmonic-2013}.

The effect of a thin layer containing wires was investigated in \cite
{Schweizer-Wiedemann-Pol-2025} with the result that a polarization
filtering interface can occur. We develop these methods further in
order to cover the critical scaling. As outlined above, the effective
system is not the one that could be expected on the basis of \cite
{Schweizer-Wiedemann-Pol-2025}. We mention that the analysis of the
effective system of \cite {Schweizer-Wiedemann-Pol-2025} is highly
non-trivial, see \cite {Delourme-Schweizer-Wiedemann-Pol-2025} for
results.

Our methods are quite classical and can well be compared with those of
\cite{MR1493040}.  We construct adequate test-functions that can be
used in the equations in order to derive effective equations. We
mention that the analysis of inclusions along a manifold was
simplified in \cite {Schweizer-Neumann-sieve-2020} for the
non-critical scaling.  The Maxwell system with non-critical wires was
analyzed in \cite {Schweizer-Wiedemann-Pol-2025} and
\cite{DelourmeHewett2020}.

Regarding the wider field of homogenization of small inclusions that
concentrate along a manifold, we mention some other works and refer
also to the references therein: \cite {Neuss-Radu-Jaeger} for
reaction-diffusion equations, \cite {MR884812} for a fluid-flow
problem, \cite {DelourmeHaddarJoly-2012} for a Helmholtz problem, and
\cite {MR1676922} for more complex geometries.

The effective law in our setting, more precisely, the interface
condition \eqref {eq:int-law-hc-intro}, is interesting due to a term
with a second derivative: $\del_1^2 j$. Such a term can produce
resonances in the effective system, but it can also be regarded as a
term that is created by a resonance: The small structure can interact
with the surrounding field and can, even though the wires are
exponentially thin, produce terms in the effective system. Regarding
effective systems that are the result of resonances, we mention \cite
{Donato-Lamacz-Schweizer-2022}, \cite {Lamacz-Schweizer-2017} and
\cite {Maurel-etal-2019} for Helmholtz-resonator resonances, we
mention again \cite {Bouchitte-Schweizer-Max-2010} for a
Maxwell-system resonance and refer to \cite {Schw-resonance-2017} for
an overview of some of these results. Finally, we refer to \cite
{MR3460101} for a resonant system in elastodynamics.

\subsection{Geometry and notations}\label{sec:Geometry}

Let us specify the geometry; we keep the setting simple so that we can
concentrate on the important interface condition. The macroscopic
domain $\Omega$ is separated by an interface $\Gamma$ into an upper
part $\Omega_+$ and a lower part $\Omega_-$:
\begin{align}\label{eq:def:Omega-Gamma}
  \Omega &\coloneqq (0,1)^2 \times (-1,1) \,, & \Gamma
  &\coloneqq (0,1)^2 \times \{0\}\,,\\\label{eq:def:Omega-pm}
  \Omega_+ &\coloneqq (0,1)^2 \times (0,1) \,,&\Omega_-
  &\coloneqq (0,1)^2 \times (-1,0) \,.
\end{align}
We denote the normal vector of $\Gamma$ pointing from $\Omega_+$ to
$\Omega_-$ by $\nu$ such that $\nu \coloneqq -e_3$.

Let $\eta >0$ be a small parameter.  We study the equations for a
sequence of sets $\Sigma_\eta \subset \Omega$, which represent
perfectly or highly conducting inclusions, i.e., the wires. The sets
$\Sigma_\eta$ concentrate at the interface $\Gamma = \{x_3 = 0\}$ in
the sense that $\Sigma_\eta \subset (0,1)^2 \times (0, \eta)$.  We are
interested in the effective interface condition along $\Gamma$.

\smallskip We denote the three dimensional unit cell by
$Y \coloneqq (0,1)^3$. The microscopic geometry is prescribed by a
sequence of compact obstacle shapes
\begin{equation*}
  \Sigma_Y^\eta \subset [0,1]^2\times (0,1) \subset \overline{Y}\,.
\end{equation*}
Scaled versions of this set form a subset of $\Omega$ and we decompose
$\Omega$ accordingly,
\begin{equation*}
  \Sigma_\eta \coloneqq
  \Omega\cap  \bigcup\limits_{k \in \Z^2}
  \eta[(k_1,k_2,0) + \Sigma_Y^\eta]\,,
  \qquad 
  \qquad 
  \Omega_\eta \coloneqq \Omega \setminus \Sigma_\eta \,.
\end{equation*}

\begin{figure}[ht]
  \centering
  \begin{subfigure}[t]{.48\textwidth}
    \centering
    \raisebox{0.4cm}{\includegraphics[width=0.8\linewidth]{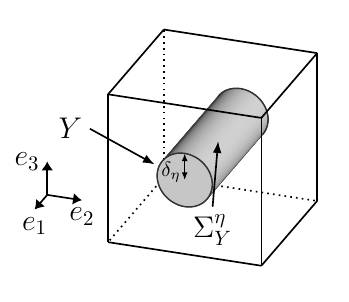}
    }
    \captionsetup{aboveskip=2pt}
    \subcaption{The set $\Sigma_Y^\eta$ representing the conductor
      in the reference cell}
    \label{fig:Y-Sigma}
  \end{subfigure}
  \hspace{\fill}
  \begin{subfigure}[t]{.48\textwidth}
    \centering
    \includegraphics[width=0.8\linewidth]{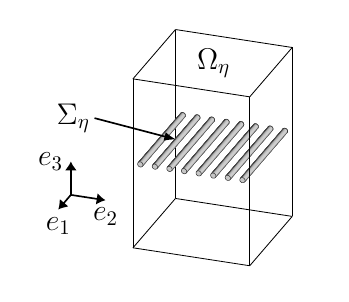}
    \captionsetup{aboveskip=4pt} \subcaption{The obstacle set
      $\Sigma_\eta$ and the domain $\Omega_\eta$}
    \label{fig:Omega_Eta}
  \end{subfigure}
  \captionsetup{aboveskip=0pt}
  \caption{Geometry of the conductor}
  \label{fig:geometries}
\end{figure}

Our analysis focusses on wire type obstacles $\Sigma^\eta_Y$ that are
aligned in the direction $e_1$ and have a small radius
$\delta_\eta>0$. We allow that the cross-section is given by a simply
connected closed set $\cC_{2D}\subset \bar B_1(0)\subset [-1,1]^2$
with a $C^{2,\gamma}$-boundary curve (we assume a high regularity in
order to simplify the arguments that are connected to extensions of
cell-solutions into the inclusion).  With the central point
$z_0 = (1/2, 1/2) \in \R^2$, we use
\begin{align}
  \label{eq:wire-description}
  \Sigma^\eta_Y \coloneqq [0,1] \times (z_0 + \delta_\eta\, \cC_{2D})
\end{align}
One may always think of circular cross-sections that are realized with
the choice $\bar B_{1}(0)\subset\R^2$. Later on, we also write
$\cC_\eta \coloneqq (0,1)\times (z_0 + \delta_\eta\, \cC_{2D})$.  For
an illustration of the cell-geometry see Figure \ref{fig:Y-Sigma}, the
resulting sets $\Sigma_\eta$ and $\Omega_\eta$ are shown in
Figure~\ref{fig:Omega_Eta}.

\vspace*{-3mm}

\captionsetup[subfigure]{aboveskip=-8pt, belowskip=8pt,
  labelfont=up, labelformat=simple}
\renewcommand\thesubfigure{(\alph{subfigure})}

For a measurable set $U \subset \R^n$, we denote the indicator
function by $\one{U}$, i.e.~${\one{U}(x) = 1}$ for $x \in U$ and
$\one{U}(x) = 0$ for $x \not \in U$. In estimates, we use generic
constants $C$ that may change from one inequality to the next; they
are always independent of the parameter $\eta$.

\subsection{Weak form of the equations}
For the analysis of Maxwell's equations we use the following function
spaces: For a domain $U\subset \R^3$, we understand the curl of an
$L^2(U, \C^3)$-function in the sense of distributions and set
\begin{align*}
  H(\rot, U)
  &\coloneqq \{ E \in L^2(U,\C^3) \mid \rot E \in L^2(U,\C^3)\} \,,\\
  H_0(\curl,U)
  &\coloneqq \left\{E\in H(\curl,U)\,\middle|\,\int_U 
    E\cdot\curl\psi
    = \int_U \curl E\cdot \psi\ \forall
    \psi\in H(\curl,U)\right\} \,,\\
  H_\loc(\rot, U)
  &\coloneqq \{ E \in L^2_\loc(U,\C^3) \mid \rot E \in L^2_\loc(U,\C^3)\}\,.
\end{align*}
For the first two spaces, we use the norm and the scalar product that
are induced by
\begin{align*}
  \|E \|^2_{H(\rot, U)} \coloneqq
  \int_U \left\{|E|^2 + |\rot E|^2\right\} \, .
\end{align*}
We use the following weak formulation of the $\eta$-problem
\eqref{eq:Strong:Maxwell:eta} for highly conducting wires: Find
$(\Ee, \He) \in L^2(\Omega, \C^3) \times L^2(\Omega, \C^3)$ such that
\begin{subequations}\label{eq:weak:Eta:E-H}
  \begin{align}\label{eq:weak:Eta:E}
    \int_{\Omega} \Ee \cdot \curl \psi
    &= \int_{\Omega} (i \omega \mu \He + f_h) \cdot \psi
    && \text{ for every } \psi \in H(\curl, \Omega)\,,
    \\\label{eq:weak:Eta:H}
    \int_{\Omega} \He \cdot \curl \phi
    &= \int_{\Omega} (-i \omega \eps_\eta \Ee + f_e) \cdot \phi
    && \text{ for every } \phi \in H_0(\curl, \Omega)\,.
  \end{align}
\end{subequations}
We note that \eqref{eq:weak:Eta:E} implies that
\eqref{eq:Strong:Maxwell:eta:1} is satisfied in the sense of
distributions (and, hence, also the regularity
$\curl \Ee \in L^2(\Omega, \C^3)$ is satisfied), \eqref{eq:weak:Eta:E} implies
also the boundary condition, $\Ee\in H_0(\curl, \Omega)$ and, hence,
\eqref {eq:Strong:Maxwell:eta:3} in a weak form.

The weak form of the $\eta$-system for perfect conductors \eqref
{eq:Maxwell:eta:pc} is given by \eqref {eq:weak:Eta:E-H} with all
integrals over $\Omega$ replaced by integrals over $\Oe$ and the
unknowns being in the function spaces with $\Omega$ replaced by $\Oe$.

\subsection{Main result}

Our main contribution is the derivation of the missing interface
relation \eqref{eq:Strong:Maxwell:limit:int} for the critical scaling.

\begin{theorem}[Effective Maxwell system with interface condition]
  \label{thm:main}
  Let the geometry be given by $\Omega$ and $\Sigma_\eta$ as described
  in Section \ref {sec:Geometry} with the critical scaling, see \eqref
  {eq:eps-eta} and \eqref {eq:delta-theta}.  For
  $f_e,f_h \in L^2(\Omega, \C^3)$ with $\div f_e \in L^2(\Omega, \C)$,
  let $(\Ee, \He)$ be a sequence of weak solutions to the Maxwell
  equations, either for perfect conductors as in \eqref
  {eq:Maxwell:eta:pc} or with highly conducting media as in \eqref
  {eq:Strong:Maxwell:eta}.  In both cases, we assume that the sequence
  $(\Ee, \He)$ satisfies a uniform $L^2$-bound, for perfect
  conductors, this demands that
  \begin{equation}
    \label{eq:uniform-L2-bound}
    \int_{\Omega_\eta} |\Ee|^2 + |\He|^2 \le C_0
  \end{equation}
  for some $C_0$, independent of $\eta$. We choose a subsequence and
  limit functions $(E,H)$ such that the weak convergence
  $(\Ee, \He) \weakto (E,H)$ holds in
  $L^2(\Omega, \C^3)\times L^2(\Omega, \C^3)$, see \eqref
  {eq:weak-limits}; in the case of perfect conductors, we consider the
  trivial extension of the solutions to all of $\Omega$.

  Then, the limit functions $(E,H)$ satisfy the following effective
  equations:
  \begin{equation}\label{eq:weak:limit:E}
    \int_{\Omega} E \cdot \curl \psi
    = \int_{\Omega} (i \omega \mu H + f_h) \cdot \psi
  \end{equation}
  for every $\psi \in H(\curl, \Omega)$. The equation for $\curl H$
  includes an additional term,
  \begin{equation}\label{eq:weak:limit:H}
    \int_{\Omega} H \cdot \curl \phi
    = \int_{\Omega} (-i \omega \eps_0 E + f_e) \cdot \phi
    - \int_\Gamma (i \omega)\, j\, \phi_1
  \end{equation}
  for every $\phi \in H_0(\curl, \Omega) \cap C^2(\bar\Omega)$.  The
  limiting flux $j$ is a distribution on $\Gamma$; with a slight abuse
  of notation we write the application of the distribution $j$ to a
  test-function $\phi_1$ as an integral as in \eqref
  {eq:weak:limit:H}, compare Remark \ref{rem:Distribution-j}. Relation
  \eqref {eq:weak:limit:H} implies that $j$ can be expressed in terms
  of the traces of $H_2$ through
  $\llbracket H_2 \rrbracket_\Gamma = i \omega\, j$.
  
  The remaining interface condition demands a relation between the
  distributions $j$ and $E_1$. We use three parameters: $\beta>0$ is
  related to the volume of the cross-section of the wires, $\kappa>0$
  is related to the conductivity, and $\alpha>0$ is determined by the
  size of the inclusions.  For perfect conductors, there holds
  \begin{equation}
    \label {eq:int-law-pc}
    \left(\omega^2 \mu\eps_0\, +\, \del_1^2\right)\, j
    = -\frac{2\eps_0}{\alpha} E_1|_\Gamma\,.
  \end{equation}
  For highly conducting wires, there holds
  \begin{equation}\label {eq:int-law-hc}
    \left(\omega^2 \mu\eps_0\, +\, \del_1^2\right)\, j
    + \frac{2i \eps_0}{\alpha\beta\kappa}\, j
    = -\frac{2\eps_0}{\alpha} E_1|_\Gamma\,.
  \end{equation}
  The distribution $j$ satisfies a Neumann condition in the sense
  that, for test-functions $\vth$ that do not vanish on
  $\{ x_1 = 0 \}$ or $\{ x_1 = 1 \}$ (but with $\del_1\vth = 0$
  there), the application of $\del_1^2 j$ to $\vth$ is given by
  $- \int_\Gamma \del_1 j\, \del_1 \vth$.
\end{theorem}

\begin{proof}[Easy parts of the proof]
  The relation \eqref {eq:weak:limit:E} is obtained by taking the
  limit $\eta\to 0$ in \eqref {eq:weak:Eta:E}, which holds in the same
  form also for perfect conductors when the functions are trivially
  extended.
  
  Relation \eqref {eq:weak:limit:H} for functions $\phi$ that vanish
  along $\Gamma$ follows directly from \eqref {eq:weak:Eta:H}.  This
  implies that, when we define a distribution $\vec j$ on $\Gamma$
  through the relation
  $i\omega \vec j = \llbracket H\times e_3 \rrbracket_\Gamma$,
  relation \eqref {eq:weak:limit:H} holds for arbitrary smooth
  functions $\phi$ with $j$ replaced by $\vec j$. We see that \eqref
  {eq:weak:limit:H} contains the additional information that
  $\vec j = j\, e_1$ for some scalar distribution $j$. For perfect
  conductors, this follows from the result
  $\llbracket H_1 \rrbracket_\Gamma=0$, for which we refer to \cite
  {Schweizer-Wiedemann-Pol-2025}. For highly conducting wires, this
  fact is shown in Lemma \ref {lem:Jj} (together with Appendix \ref
  {app.limit-J}).

  It remains to show \eqref {eq:int-law-pc} and \eqref
  {eq:int-law-hc}. This is our goal in the next three sections.
\end{proof}

\begin{remark}\label{rem:Distribution-j}
  The tangential trace of a function in $H(\curl)$ is an element in
  $H^{-1/2}(\operatorname{Div})$ on the boundary of the domain, but,
  in general, not an $L^2$-function. When we write
  $\llbracket H_2 \rrbracket_\Gamma = i \omega\, j$, this must be
  interpreted in the sense of distributions,
  \begin{equation*}
    \la i\omega\, j, \phi_1 \ra
    = \la \llbracket H_2 \rrbracket_\Gamma, \phi_1 \ra
    \coloneqq
    \int_{\Omega\setminus \Gamma} \left\{\curl H \cdot (\phi_1, 0, 0)
      - H \cdot \curl (\phi_1,0,0) \right\}
  \end{equation*}
  for every $\phi_1\in C^2(\bar\Omega)$ such that
  $(\phi_1,0,0)\in H_0(\curl,\Omega)$.  Our (sloppy) notation is to
  write $i\omega \int_\Gamma j \phi_1$ for the left-hand side.  As
  usual, derivatives of $j$ such as $\del_1^2 j$ are defined by
  letting the derivative act on the test-function.  Similarly,
  $E_1|_{\Gamma}$ in \eqref {eq:int-law-pc} and \eqref{eq:int-law-hc}
  is the distribution defined by
  \begin{equation*}
    \la E_1|_\Gamma , \phi_2 \ra
    \coloneqq -\frac{1}{2}\int_\Omega \sign(x_3)
    \left\{ \curl E \cdot (0, \phi_2,0) - E \cdot \curl (0, \phi_2,0)\right\}\,.
  \end{equation*}
  Moreover, the interface condition
  \eqref{eq:Strong:Maxwell:limit:int-triv} has to be understood in the
  distributional sense, where $\llbracket H_1\rrbracket_\Gamma$ is
  given by
  \begin{equation*}
    \label{eq:def:H1}
    \la \llbracket H_1\rrbracket_\Gamma, \phi_2 \ra \coloneqq
    - \int_\Omega \left\{\curl H \cdot (0,\phi_2, 0)
      - H \cdot \curl (0,\phi_2,0) \right\}\,.
  \end{equation*}

\end{remark}

\section{Cell functions for wires}

We follow the typical route of homogenization in critical scaling: We
define cell-functions, the underlying idea is that these cell
functions should look like solutions in the vicinity of the
inclusions. Using the cell functions as test-functions in the
equations, we find, in the limit $\eta\to 0$, the homogenized system.

\subsection{Singular cell-functions for perfect conductors}

Shape functions are defined on cylinders $Z = (0,1)^2\times \R$. We
regard the central domain as the periodicity cell,
$Y = (0,1)^3\subset Z$.  The wire defines a cylindrical inclusion
$\cC_\eta \subset Z$. In the simplest case of circular cross-sections,
the wires can be chosen as given by
$\cC_{2D} = \bar B_{1}(0) \subset \R^2$ and
$\cC_\eta = (0,1)\times \bar B_{\delta_\eta}(z_0) \subset Z$ with
$z_0 = (1/2, 1/2) \in \R^2$, see \eqref {eq:wire-description}.

In our first definition of shape-functions, we define test-functions
outside the wires. The idea is that the functions mimic the expected
behavior of solutions. Only the shape-functions of Definition \ref
{def:shapeFunction} and their less-singular counterparts of Definition
\ref {def:less-singular} are required for the analysis of perfect
conductors.

\begin{definition}[Shape-functions outside the
  wires]\label{def:shapeFunction}
  We define the first shape-function $U^\eta \colon Z\to \C^3$ as a
  function that points in direction $e_1$, hence
  $U^\eta_2 = U^\eta_3 = 0$.  The first component
  $U^\eta_1 = U^\eta_1(y) = U^\eta_1(y_2,y_3)$ is defined as the
  unique solution of the problem
  \begin{align*}
    \Delta U_1^\eta &= 0 \text{ on } Z\setminus \cC_\eta\,,\\
    U_1^\eta(y) &= 0\text{ on } \del \cC_\eta\,,\\
    \nabla U_1^\eta(y) - \eta e_3\, \sign(y_3)
                    &  \text{ is of class  } L^2(Z\setminus \cC_\eta)\,,\\
    U_1^\eta(y) &\text{ periodic on } \{y_2=0\}\text { and } \{y_2=1\}\,.
  \end{align*}
  The above problem has a unique solution, see Appendix \ref {app.UV}.
  The solution has the property that $U_1^\eta(y) - \eta |y_3|$ is
  bounded.  The solution has, for some real numbers
  $\alpha_\eta^\pm \in \R$, the asymptotic behavior
  $U_1^\eta(y) - \eta |y_3| \to \alpha_\eta^\pm$ for
  $\pm y_3 \to \infty$.  The shifts satisfy
  $\alpha_\eta^\pm \to \alpha > 0$, the parameter $\alpha > 0$ is an
  important quantity in the effective equations.
  
  The second shape-function
  $V^\eta \colon Z \setminus \cC_\eta \to \C^3$ is defined in terms of
  the curl of the first shape-function:
  \begin{equation}
    V^\eta \coloneqq \eta^{-1} \curl(U^\eta)\,.
    \label{eq:V-eta}  
  \end{equation}
\end{definition}

{\bf Remarks regarding the qualitative behavior.}  The singular
behavior of $U_1^\eta$ near $z_0$ is characterized by the formula
\begin{equation*}
  U_1^\eta(y) \sim \eta \pi^{-1}\, \log(r/\delta_\eta)\, \qquad
  \text{ for } \ r \coloneqq |(y_2,y_3) -z_0|\,. 
\end{equation*}
Indeed, when we consider radii $r$ of order $1$, this function is of
order $1$: Using the critical scaling of $\delta_\eta$ in
\eqref{eq:delta-theta}, we can calculate
$\eta \pi^{-1} \log(1/\delta_\eta) = \eta\pi^{-1} \log( e^{\alpha\pi/\eta} ) =
\alpha$.

The second shape function has the following features: As a curl, it
has vanishing divergence. Just as $U^\eta$, it depends only on
$(y_2, y_3)$ and it is periodic in $y_2$. It points everywhere in the
$(y_2, y_3)$-plane, there holds $V^\eta_1 = V^\eta\cdot e_1 =
0$. Since $U^\eta$ also has vanishing divergence, there holds
$\curl V^\eta = \curl \curl U^\eta = -\Delta U^\eta = 0$ in
$Z\setminus \cC_\eta$.

Away from $\cC_\eta$, derivatives of the function $U^\eta$ are
small. More precisely, on the set $Z\setminus B_{1/2}(z_0)$,
derivatives are at most of the order
$\del_r U_1^\eta(y) = \del_r [\eta\, \log(r/\delta_\eta)] \sim \eta$.
Therefore $V^\eta$ is bounded on the set $Z\setminus B_{1/2}(z_0)$,
independent of $\eta$. We note that the sequence $V^\eta$ is bounded
in $L^1$ on bounded sets, but it is not bounded in $L^2$ on bounded
sets.

In the vicinity of the inclusion $\cC_\eta$, the shape function is of
the order $r^{-1}$. More precisely, for circular cross-sections, the
shape is of the qualitative form
\begin{equation*}
  V^\eta(z_0 + y) \sim \frac{(0, y_3, -y_2)}{|y_2|^2+|y_3|^2}\,.
\end{equation*}
Far away from the obstacle, i.e., in the limit $|y_3| \to \infty$, the
field stabilizes: There holds
\begin{equation}
  \label{eq:V-stabilization}
  V^\eta(y) \to \pm e_2\quad \text{ as } |y_3|\to \infty\,.
\end{equation}
We verify also this property in Appendix \ref {app.UV}.

\smallskip In the original problem, posed in the domain $\Omega_\eta$,
we use test-functions $\Phi^\eta$ and $\Psi^\eta$ that are scaled
versions of $U^\eta$ and $V^\eta$.

\begin{definition}[Scaled cell-functions]
  \label{def.sing-cell-fcts}
  We define cell-functions $\Phi^\eta, \Psi^\eta \colon \Oe \to \R^3$ as
  \begin{equation}
    \label{eq:cell-functions-PP}
    \Phi^\eta(x) \coloneqq U^\eta(x/\eta)\ \text{ and }\ 
    \Psi^\eta(x) \coloneqq V^\eta(x/\eta)\ \text{ with }\ 
    \curl\Phi^\eta = \Psi^\eta\,.
  \end{equation}
\end{definition}

The following properties are satisfied, for proofs see Appendix \ref
{app.UV}: The function $\Phi^\eta$ is uniformly bounded, it points
everywhere in $e_1$-direction and has a vanishing divergence, and its
trivial extension to $\Omega$ converges strongly in
$L^2(\Omega, \C^3)$ to the function
$\Phi^0(x) \coloneqq (\alpha + |x_3|)\, e_1$.

The function $\Psi^\eta$ is everywhere orthogonal to $e_1$. The
sequence $\Psi^\eta$ is bounded in $L^2(\Oe, \C^3)$.  Its trivial
extension to $\Omega$ converges weakly to
$\sign(x_3) e_2 = \curl \Phi^0(x)$. It has a vanishing divergence and
a vanishing curl in $\Oe$.

\subsection{Extended cell-functions}

For the analysis of highly conducting wires, we extend the shape
functions into the set $\cC_\eta$.  The function $U^\eta$ vanishes
along $\del \cC_\eta$; we extend it trivially and set $U^\eta(y) = 0$
for $y\in \cC_\eta$.

We next extend $V^\eta(y)$ inside the wires in such a way that
tangential components of $V^\eta$ do not jump along $\del\cC_\eta$.
In the case that the wires have a circular cross-section, we can
essentially extend $V^\eta$ as a function that describes a rigid
rotation.  For $y = (y_2,y_3)$ with $|y| = \delta_\eta$, we expect
\begin{equation*}
  V^\eta(z_0 + y) \approx \frac{1}{\pi\delta_\eta^2} (0,y_3, -y_2)\,.
\end{equation*}
Loosely speaking, our aim is to extend $V^\eta$ with the same formula
for $|y| < \delta_\eta$. Accordingly, inside the wire, we expect
$\curl V^\eta \approx -2/(\pi\delta_\eta^2)\, e_1$.

For the rigorous construction, we use the solution of a problem in the
obstacle $\cC_\eta$. As a preparation, we note that the function
$U_1^\eta$ satisfies
\begin{equation*}
  \int_{ (0,1)^2 \times \{y_3\} }
  \del_3 U_1^\eta = \pm\eta\text{ for }
  y_3\to\pm \infty\,.
\end{equation*}
The Theorem of Gauß implies that the normal derivative along the
obstacle satisfies
\begin{equation*}
  g \coloneqq \del_\nu U_1^\eta|_{\del\cC_\eta}\,,\quad
  \int_{\del\cC_\eta} g = 2\eta\,,
\end{equation*}
where we use the convention that the normal vector $\nu$ on
$\del\cC_\eta$ points to the outside of the obstacle. We conclude that
the function $g$ has typical values of the order
$\eta/|\del\cC_\eta| \sim \eta/\delta_\eta$.

We define $G^\eta \colon \cC_\eta\to \R$ as the solution of the problem
\begin{align*}
  \Delta G^\eta &= 2\eta/|\cC_\eta|\quad \text{ in } \cC_\eta\,,\\
  \del_\nu G^\eta &= g\quad  \text{ on } \del\cC_\eta\,,
\end{align*}
with the normalization that $G^\eta$ has vanishing average. Here, we
have respected the integrability condition that is dictated by
$\int_{\cC_\eta} \Delta G^\eta = \int_{\del\cC_\eta} g$. With this
construction, we can now define $V^\eta$ inside the obstacle:
\begin{equation}
  V^\eta \colon \cC_\eta\to \R^3\,,
  \qquad V^\eta \coloneqq \eta^{-1} \curl((G^\eta,0,0))\,.
  \label{eq:V-eta-cC}
\end{equation}
The definition ensures that tangential components of $V^\eta$ have no
jump across $\del\cC_\eta$, cf. \eqref{eq:V-eta} for the definition of
$V^\eta$ outside the obstacles.  The values are bounded in the sense
that $|V^\eta| \le C/\delta_\eta$ holds pointwise with $C$ independent
of $\eta$, see Appendix \ref {app.UV}. The curl is given as
$\curl V^\eta = - \eta^{-1} \Delta G^\eta\, e_1 = - 2/|\cC_\eta|\,
e_1$.

We define the extended cell-functions as in Definition \ref
{def.sing-cell-fcts}, $\Phi^\eta, \Psi^\eta \colon \Omega \to \R^3$ are
defined as
\begin{equation}
  \label{eq:cell-functions-PP-extended}
  \Phi^\eta(x) \coloneqq U^\eta(x/\eta)\ \text{ and }\ 
  \Psi^\eta(x) \coloneqq V^\eta(x/\eta)\,.
\end{equation}
The sequence of extended cell-functions $\Psi^\eta$ is bounded in
$L^2(\Omega,\C^3)$.  Because of \eqref {eq:V-stabilization}, we can
determine its weak limit in that space,
\begin{equation}
  \label{eq:Psi-eta-weakly}
  \Psi^\eta \weakto \sign(x_3) e_2\,.
\end{equation}
The relation $\curl \Phi^\eta = \Psi^\eta$ holds outside the wires,
{\em but not} in the wires.  In $\Sigma_\eta$ holds a uniform estimate
$|\Psi^\eta| \le C/\delta_\eta$. The curl is constant in the wire,
namely 
\begin{equation}\label{eq:Curl-Psi}
  \curl \Psi^\eta = - 2/|\cC_\eta| \eta^{-1} \, e_1  \one{\Sigma_\eta}\,.
\end{equation}

\subsection{Less singular test-functions}
We define $\tilde U^\eta(y)$ as a more regular version of $U^\eta(y)$:
Vanishing in $B_{1/2}(z_0)$ and essentially equal to
$U^\eta(y) - \alpha$ outside $Y$. For instance, one may choose
$\tilde U^\eta(y_1,y_2,y_3) = e_1 \left( \eta (y_3 -1) \,\one{\{y_3
    \geq 1\}} -\eta y_3 \,\one{\{y_3 \leq 0\}}\right)$.  Note that
$\tilde{U}^\eta_1$ is not harmonic (in contrast to $U^\eta_1$).  The
construction of the function $V$ is as before, we set
$\tilde V^\eta \coloneqq \eta^{-1} \curl(\tilde U^\eta)$.

\begin{definition}[Less singular cell-functions]
  \label{def:less-singular}
  We define the cell-function $\tilde\Phi^\eta \colon \Omega \to \R^3$ as
  $\tilde\Phi^\eta(x) \coloneqq \tilde U^\eta(x/\eta)$. It vanishes in the
  wires and converges strongly in $L^2(\Omega)$ to $|x_3|\, e_1$.

  The second cell-function is $\tilde\Psi^\eta \colon \Omega \to \R^3$ as
  $\tilde\Psi^\eta(x) \coloneqq \tilde V^\eta(x/\eta)$. It is bounded in
  $L^\infty(\Omega)$ and converges strongly in $L^2(\Omega,\C^3)$ to
  $\sign(x_3) e_2$.
\end{definition}

The difference $W^\eta(y) \coloneqq U_1^\eta(y) - \tilde U_1^\eta(y)$ and the
rescaled version $w^\eta(x) \coloneqq W^\eta(x/\eta)$ allow us to rewrite one
critical contribution in terms of a gradient. We use the symbol
$(0,b,c)^\perp = (0,-c,b)$ for $b,c\in \C$. For $x\in \Oe$ there holds
\begin{align}
	\begin{aligned}
  \Psi^\eta(x) - \tilde\Psi^\eta(x)
  &= \curl \Phi^\eta(x) - \curl \tilde\Phi^\eta(x)
    = - \frac1\eta \left[\nabla_y U^\eta_1(x/\eta)
    - \nabla_y \tilde U^\eta_1(x/\eta) \right]^\perp\\
  &= - \frac1\eta \left[\nabla_y W^\eta(x/\eta)\right]^\perp
    = - \frac1\eta\ e_1\times \nabla_y W^\eta(x/\eta)
    = - e_1\times \nabla_x w^\eta(x)\,.
    \label{eq:Psi-diff-nabla-W}
\end{aligned}
\end{align}
The function $w^\eta$ converges strongly to the constant function
$\alpha$ in $L^2(\Omega, \C)$. It vanishes in $\Sigma_\eta$ and is
independent of $y_1$ and, thus, $\partial_1 w^\eta=0$.  We call
$\tilde{\Phi}^\eta$ less singular cell function since, in contrast to
the function $\Phi^\eta$, it converges strongly in $H(\curl, \Omega)$.

\section{The limit problem for perfect conductors}

In both cases, for perfect conductors and for highly conducting wires,
we want to interpret the limit of the second equation,
\eqref{eq:Maxwell:eta:pc:2} or \eqref{eq:Strong:Maxwell:eta:2}, as
$\curl H = -i\omega \eps_0 E + f_e - i \omega j e_1
\cH^2\lfloor_\Gamma$; we refer once more to Remark
\ref{rem:Distribution-j} regarding the interpretation of $j$.  This
will indeed be our result in \eqref {eq:weak:limit:H}. The first
component of this relation contains, on the left-hand side, the
singular part of $-\del_3 H_2$. Comparing the singular parts of both
sides, we therefore find
$- \llbracket H_2 \rrbracket_\Gamma = - i \omega j$.  In the case of
perfect conductors, we therefore {\em define} the current $j$ as
\begin{equation}
  \label{eq:j-perf-cond}
  j \coloneqq \frac{1}{i \omega} \llbracket H_2 \rrbracket_\Gamma\,.
\end{equation}
We note that there is a more natural definition of $j$ in the case of
highly conducting wires, namely as the limit of $\eps_\eta E^\eta_1$
in the wires. As it turns out, this other definition leads to the same
quantity as above, but we cannot use the other definition in the case
of perfect conductors. In any case, this observation explains why we
include the factor $i \omega$ in \eqref {eq:j-perf-cond}.

Let us investigate the quantity $j$ further. We consider the third
component of $\curl H = - i \omega \eps_0 E + f_e$, which is, wherever
$f_e$ vanishes, $\del_1 H_2 - \del_2 H_1 = - i \omega \eps_0 E_3$ in
$\Omega\setminus \Gamma$. We evaluate the jumps of both sides. Since
the wires are disconnected in direction $e_2$, there holds
$\llbracket H_1 \rrbracket_\Gamma = 0$ and hence also
$\llbracket \del_2 H_1 \rrbracket_\Gamma = 0$. Assuming that $f_e$
vanishes in the vicinity of $\Gamma$, evaluating the jumps, we find,
in the sense of distributions,
\begin{equation}
  \label{eq:div-E-pc}
  - i\omega \eps_0 \llbracket E_3 \rrbracket_\Gamma
  = \llbracket \del_1 H_2 \rrbracket_\Gamma
  = \del_1 \llbracket H_2 \rrbracket_\Gamma
  = i\omega \del_1 j\,.
\end{equation}

\smallskip The homogenization result is based on the following
procedure: For an arbitrary function $\vth\in C_c^\infty(\Omega)$, we
use
\begin{align*}
  \phi(x) &\coloneqq \vth(x) \Phi^\eta(x)\qquad
            \text{ in equation \eqref {eq:Maxwell:eta:pc:2}
            for } \curl(H^\eta)\,,\\
  \psi(x) &\coloneqq \vth(x) \Psi^\eta(x)\qquad
            \text{ in equation \eqref
            {eq:Maxwell:eta:pc:1} for } \curl(E^\eta)\,.
\end{align*}
We note that we can also allow that $\vth$ does not vanish on
$\{ x_1 = 0 \}$ or $\{ x_1 = 1 \}$ (but we demand for the proofs that
$\del_1\vth = 0$ is satisfied there). Indeed, no boundary conditions
are imposed in \eqref{eq:weak:Eta:E}. In \eqref{eq:weak:Eta:H}, the
test-function $\phi$ must have vanishing tangential components; this
is guaranteed since $\Phi^\eta$ points in $e_1$-direction.

To simplify the calculations, we always assume that $f_e$ and $f_h$
have a support that does not intersect $\Gamma$. Under this
assumption, we can choose $\vth$ such that its support is disjoint
from the supports of $f_e$ and $f_h$. This does not restrict the
generality of the argument, we simply avoid the terms connected to
$f_e$ and $f_h$ in the subsequent calculations.

Using $\phi = \vth \Phi^\eta$ in \eqref {eq:Maxwell:eta:pc:2}, an
application of the product rule gives
\begin{equation}
  \label{eq:testedH-pc}
  - \int_\Oe i \omega \eps_0 E^\eta\cdot \vth \Phi^\eta
  = \int_\Oe \vth\, \He \cdot \curl \Phi^\eta 
  - \int_\Oe \He \cdot \Phi^\eta \times \nabla\vth\,.
\end{equation}
Using $\psi = \vth \Psi^\eta$ in \eqref {eq:Maxwell:eta:pc:1}, we can
exploit $\curl\Psi^\eta = 0$ and find
\begin{equation}
  \label{eq:testedE-pc}
  \int_\Oe i \omega \mu H^\eta\cdot \vth \Psi^\eta
  = - \int_\Oe \Ee \cdot \Psi^\eta \times \nabla\vth\,.
\end{equation}
We constructed test-functions with $\curl\Phi^\eta = \Psi^\eta$. This
property implies that the second integral in \eqref {eq:testedH-pc}
coincides with the left-hand side of \eqref {eq:testedE-pc} up to the
factor $i\omega\mu$. We can therefore combine the two equations \eqref
{eq:testedH-pc}--\eqref {eq:testedE-pc} to
\begin{equation}
  \label{eq:combined-EH-p}
  - \int_\Oe i \omega \eps_0 E^\eta\cdot \vth \Phi^\eta
  + \int_\Oe \He \cdot \Phi^\eta \times \nabla\vth
  = \frac{i}{\omega\mu} \int_\Oe \Ee \cdot \Psi^\eta \times \nabla\vth\,.
\end{equation}
Denoting the integrals by $I_{1}^\eta$ to $I_{3}^\eta$, we write this
relation as $I_{1}^\eta + I_{2}^\eta = I_{3}^\eta$. It remains to
determine the limits of the integrals as $\eta\to 0$. This will
provide \eqref {eq:int-law-pc}.

\begin{lemma}[The limits of $I_{1}^\eta$ and $I_{2}^\eta$]
  \label{lem:easy-1p}
  There holds, in the limit $\eta\to 0$,
  \begin{equation}
    \label{eq:easy-1p}
    I_{1}^\eta \coloneqq - \int_\Oe i \omega \eps_0 E^\eta\cdot \vth \Phi^\eta
    \to - \int_\Omega i \omega \eps_0 E\cdot \vth \,(\alpha + |x_3|)\, e_1\,,
  \end{equation}
  and
  \begin{equation}
    I_{2}^\eta = \int_\Oe \He \cdot \Phi^\eta \times \nabla\vth
    \to \int_\Omega H \cdot (\alpha + |x_3|)\, e_1 \times \nabla\vth\,.
    \label{eq:easy-2p}
  \end{equation}
\end{lemma}

\begin{proof}
  Both limits follow easily from the strong $L^2$-convergence
  $\Phi^\eta(x) \to \Phi^0(x) = (\alpha + |x_3|)\, e_1$. The functions
  $E^\eta$ and $H^\eta$ (or, more precisely, their trivial extensions
  to all of $\Omega$) converge weakly in $L^2(\Omega, \C^3)$. This
  allows us to conclude the convergence of the product.
\end{proof}

\begin{lemma}[The limit of $I_{3}^\eta$]
  \label{lem:easy-3p}
  There holds, in the limit $\eta\to 0$,
  \begin{equation}
    \label{eq:easy-3p}
    (-i\omega\mu)\, I_{3}^\eta = \int_\Oe \Ee \cdot \Psi^\eta \times \nabla\vth
    \to \int_\Omega  E \cdot \sign(x_3) e_2 \times \nabla\vth
    - \frac\alpha{\eps_0} \int_\Gamma j\, \del_1^2 \vth\,.
  \end{equation}
\end{lemma}

\begin{proof}
  We recall that $\Psi^\eta$ converges weakly to $\sign(x_3) e_2$, see
  \eqref {eq:Psi-eta-weakly}, but this weak convergence is not
  sufficient for the limit procedure. We can proceed when we replace
  $\Psi^\eta$ by $\tilde \Psi^\eta$ -- the limit can then be performed
  by the strong $L^2(\Omega, \C^3)$-convergence of $\tilde \Psi^\eta$,
  \begin{equation}
    \label{eq:easy-3p-89654}
    \int_\Oe \Ee \cdot \tilde\Psi^\eta \times \nabla\vth
    \to \int_\Omega E \cdot \sign(x_3) e_2 \times \nabla\vth
    = - \int_\Omega E \cdot \sign(x_3) \curl (e_2 \vth) \,.
  \end{equation}
  This yields the first term on the right-hand side of \eqref
  {eq:easy-3p}.
  
  It remains to study the difference of the left-hand sides in \eqref
  {eq:easy-3p} and \eqref {eq:easy-3p-89654}.  We use the function
  $w^\eta$ with the property \eqref {eq:Psi-diff-nabla-W} to calculate
  \begin{align*}
    \int_\Oe \Ee \cdot \left(\Psi^\eta - \tilde\Psi^\eta\right) \times \nabla\vth
    &= -\int_\Oe \Ee \cdot (e_1\times \nabla w^\eta) \times \nabla\vth\,.
  \end{align*}
  We use the BAC-CAB rule in the form
  $(b \times c)\times a = c(a\cdot b)-b(a\cdot c)$ with the vectors
  $b = e_1$, $c = \nabla w^\eta$ and $a = \nabla\vth$. This allows us to
  re-write the above result as
  \begin{equation}
    \label{eq:BAC-CAB-pc}
    \int_\Oe \Ee \cdot \left(\Psi^\eta - \tilde\Psi^\eta\right) \times \nabla\vth
    = - \int_\Oe \Ee\cdot  \nabla w^\eta\  \del_1\vth
    + \int_\Oe \Ee_1\  (\nabla w^\eta\cdot \nabla\vth)\,.
  \end{equation}
  We start with an analysis of the first integral on the right. We
  integrate by parts, using that $w^\eta$ vanishes along the boundary
  $\del\Sigma_\eta$ and $\del_1\vth$ on the macroscopic
  $x_1$-boundaries. The divergence of $E^\eta$ satisfies
  $\div E^\eta = (i\omega \eps_0)^{-1}\div(f_e)$; as the support of
  $f_e$ and $\vth$ are disjoint, the corresponding contribution
  vanishes. The first integral is therefore
  \begin{align*}
    -\int_\Oe \Ee\cdot  \nabla w^\eta\  \del_1\vth
    &= \int_\Oe \Ee\cdot  \nabla \del_1\vth\ w^\eta
    \to \alpha\int_\Omega E\cdot \nabla \del_1\vth\\
    &= - \alpha\int_\Gamma \llbracket E_3  \rrbracket_\Gamma\,  \del_1\vth 
    = - \frac\alpha{\eps_0} \int_\Gamma j\,  \del_1^2\vth \,.
  \end{align*} 
  In the limit, we used that $w^\eta$ converges strongly to the
  constant function $\alpha$, see the comment after the calculation
  \eqref {eq:Psi-diff-nabla-W}.  In the last step, we used
  $\eps_0 \llbracket E_3 \rrbracket_\Gamma\, \cH^2\lfloor_\Gamma =
  -\del_1 j\, \cH^2\lfloor_\Gamma$ in the sense of distributions, see
  \eqref {eq:div-E-pc}. We conclude that the first integral on the
  right-hand side of \eqref {eq:BAC-CAB-pc} provides the desired
  surface integral in \eqref {eq:easy-3p}.
  
  The second integral on the right-hand side of \eqref {eq:BAC-CAB-pc}
  is vanishing in the limit; this follows with Lemma \ref
  {lem:int-scalar-prod} below, applied with $\psi = \nabla\vth$.
\end{proof}

We note that the last term in the above calculation was actually
obtained in the form
$(\alpha/\eps_0) \int_\Gamma \del_1 j\, \del_1\vth$. When $\vth$ is
not vanishing at the $e_1$-boundaries, this shows that $j$ satisfies a
Neumann condition.

\begin{lemma}[Limit of an error term]
  \label {lem:int-scalar-prod} 
  Let $E^\eta$ be either the solution for perfect or highly conducting
  obstacles. Then, for every smooth test-function
  $\psi\colon\Omega\to \C^3$ there holds
  \begin{equation}
    \label{eq:BAC-CAB-pc-proof}
    \int_\Oe \Ee_1\, \nabla w^\eta \cdot \psi \to 0\,.
  \end{equation}
\end{lemma}

\begin{proof}
  We define two vector fields $\sigma^\eta_{(2)}$ and
  $\sigma^\eta_{(3)}$ and compute their curl, exploiting
  $\partial_1 w^\eta = 0$:
  \begin{equation*}
    \sigma^\eta_{(j)}(x) \coloneqq (w^\eta(x) - \alpha)\, (e_1 \times e_j)\,,
    \quad \curl \sigma^\eta_{(j)}(x) = \del_j w^\eta(x)\, e_1
    \qquad \text{for }j \in \{2,3\}\,.
  \end{equation*}
  Using once more that $\partial_1 w^\eta$ vanishes, we find
  $\nabla w^\eta \cdot \psi = e_1\cdot (\curl \sigma^\eta_{(2)}\, \psi_2 +
  \curl \sigma^\eta_{(3)}\, \psi_3)$.  In the case of perfectly
  conducting obstacles, let $E^\eta$ be extended to $\Omega$ by $0$.
  Since $w^\eta$ and therefore $\nabla w^\eta$ vanish in the
  obstacles $\Sigma_\eta$, we get:
  \begin{equation*}
    \int_\Oe \Ee_1\, \nabla w^\eta\cdot  \psi     
    = \int_\Omega \Ee_1\, \nabla w^\eta\cdot  \psi     
    =
    \sum_{j\in \{2,3\}}\int_\Omega \Ee\cdot \curl
    \sigma^\eta_{(j)}\, \psi_j\,.
  \end{equation*}
  We use integration by parts, the fact that $E^\eta$ and
  $\curl E^\eta = i\omega \mu H^\eta +f_h$ are bounded in $L^2(\Omega, \C^3)$
  and the strong convergence of $w^\eta \to \alpha$ in $L^2(\Omega)$
  to
  compute
  \begin{equation*}
    \int_\Omega \Ee\cdot \curl
    \sigma^\eta_{(j)}\, \psi_j 
    =
    \int_\Omega \curl\Ee\cdot  \sigma^\eta_{(j)} \, \psi_j 
    + 
    \int_\Omega \nabla\psi_j\times \Ee\cdot \sigma^\eta_{(j)} 
    \to 0\,.
  \end{equation*}
  Combining the last two computations yields \eqref
  {eq:BAC-CAB-pc-proof}.
\end{proof}

The above lemmas provide with \eqref {eq:easy-1p}--\eqref {eq:easy-3p}
the limits of all the terms in $I_{1}^\eta + I_{2}^\eta = I_{3}^\eta$,
see \eqref {eq:combined-EH-p}. Inserting \eqref {eq:easy-1p}--\eqref
{eq:easy-3p}, we obtain
\begin{equation}
  \label{eq:res-perf-cond}
  \begin{split}
    &- \int_\Omega i \omega \eps_0 E\cdot \vth \,(\alpha + |x_3|)\, e_1
    + \int_\Omega H \cdot (\alpha + |x_3|)\, e_1 \times \nabla\vth\\
    &\quad = \frac{i}{\omega \mu} \int_\Omega E \cdot \sign(x_3) e_2 \times \nabla\vth
    - \frac{i \alpha}{\omega \mu \eps_0} \int_\Gamma j\, \del_1^2 \vth\,.
  \end{split}
\end{equation}

The bulk terms contain some well-known contributions; they are
generated by
\begin{align*}
  &\text{ \phantom{and} testing }\quad -i \omega\eps_0 E = \curl H \quad
  \text{ with }\quad \vth\, (\alpha + |x_3|)\, e_1\,,\\
  &\text{ and testing }\quad \curl E = i \omega \mu H\quad
  \text{ with }\quad \frac{i}{\omega \mu} \vth\, \sign(x_3) e_2\,.
\end{align*}
In this process, the integral $\int_\Omega H\cdot e_2 \sign(x_3) \vth$
appears twice and cancels.

Integrations by parts in \eqref {eq:res-perf-cond} introduce also
boundary terms. Comparing the boundary integrals, we find
\begin{equation}
  \label{eq:res-perf-cond-Gamma}
  \int_\Gamma \llbracket H_2 \rrbracket_\Gamma\,\alpha\vth
  = -\frac{i}{\omega \mu} \int_\Gamma E\cdot e_1\, 2 \vth
  - \frac{i \alpha}{\omega \mu \eps_0} \int_\Gamma j\, \del_1^2 \vth\,.
\end{equation}
We multiply with $i \omega \mu\eps_0/\alpha$, and use, in the first
term, relation \eqref {eq:j-perf-cond} in the form
$i\omega \llbracket H_2 \rrbracket_\Gamma = -\omega^2 j$.  Exploiting
that $\vth$ is an arbitrary function, we arrive at
\begin{equation}
  \label{eq:bdry-cond-pc}
  -\omega^2 \mu\eps_0\, j - \del_1^2\, j
  = \frac{2 \eps_0}{\alpha}\, E_1|_\Gamma\,.
\end{equation}
With this relation, we have derived \eqref {eq:int-law-pc}.

\section{The limit problem for highly conducting materials}

The case of highly conducting materials allows a more natural
definition of a flux quantity that was also used, e.g., in \cite
{Bouchitte-Schweizer-Max-2010}. The natural way to define a limiting
flux distribution $\vec J$ on $\Gamma$ is via the limit relation
\begin{equation}
  \label{eq:J-def}
  \int_\Omega \eps_\eta E^\eta\cdot \fhi
  \to \int_\Omega \eps_0 E\cdot \fhi
  + \int_\Gamma \vec J\cdot\fhi
\end{equation}
for arbitrary $\fhi\in C_c^\infty(\Omega,\R^3)$ (nonvanishing values
of $\fhi_1$ on $\{ x_1 = 0 \}$ or $\{ x_1 = 1 \}$ are also permitted).
The existence of a limiting measure $\vec J$ can be concluded from
$L^1$-estimates for $\eps_\eta E^\eta$. Furthermore, it is easy to see
that it points in direction $e_1$: There holds $\vec J = J\, e_1$ for
a complex-valued measure $J$ on $\Gamma$. Both facts follow with
standard calculations, see Appendix \ref {app.limit-J}.

We next observe that our definition of a flux with jumps of $H$,
namely through
$i\omega \vec j = \llbracket H\times e_3 \rrbracket_\Gamma$ (which is
the only possible definition for perfect conductors), coincides with
the above definition.

\begin{lemma}[The flux quantities coincide]
  \label{lem:Jj}
  There holds $\vec J = \vec j$.
\end{lemma}

\begin{proof} We use \eqref {eq:weak:Eta:H} with an arbitrary smooth
  test-function $\fhi\in C_c^\infty(\Omega,\C^3)$ to calculate
  \begin{align*}
    &\int_{\Omega} (-i \omega \eps_\eta \Ee + f_e) \cdot \fhi
      =  \int_{\Omega} \He \cdot \curl \fhi
      \to \int_{\Omega} H \cdot \curl \fhi\\
    &\quad = \int_{\Omega\setminus \Gamma}  \curl H \cdot \fhi
      - \int_{\Gamma} \llbracket H\times e_3 \rrbracket_\Gamma\cdot \fhi\\
    &\quad = \int_{\Omega} (-i \omega \eps_0 E + f_e) \cdot \fhi
    -i \omega \int_{\Gamma} \vec j \cdot \fhi\,.
  \end{align*}
  Upon subtracting the term containing $f_e$ and dividing by
  $-i\omega$, this is \eqref {eq:J-def}.
\end{proof}

We recall that the function $\Phi^\eta$ is defined in such a way that
it vanishes in the wires.  We furthermore recall that we extended the
function $\Psi^\eta$ with continuous tangential components into the
wires. To simplify the calculations, we assume once more that the
support of $\vth$ does not intersect the support of $f_e$ or $f_h$.

In equation \eqref {eq:Strong:Maxwell:eta:2}, we use the test-function
$\phi= \vth \Phi^\eta$, which provides
\begin{equation}
  \label{eq:testedH-hc}
  - \int_\Omega i \omega  \eps_\eta E^\eta\cdot \vth \Phi^\eta
  = \int_\Omega \vth\, \He \cdot \curl \Phi^\eta 
  - \int_\Omega \He \cdot \Phi^\eta \times \nabla\vth\,.
\end{equation}
In \eqref {eq:Strong:Maxwell:eta:1}, we use $\psi = \vth \Psi^\eta$ as
test-function and exploit $\curl\Psi^\eta = 0$ in $\Oe$. We find
\begin{equation}
  \label{eq:testedE-hc}
  \int_\Omega i \omega \mu H^\eta\cdot \vth \Psi^\eta
  =  \int_{\Sigma_\eta} \vth\, \Ee \cdot \curl \Psi^\eta
  - \int_\Omega \Ee \cdot \Psi^\eta \times \nabla\vth\,.
\end{equation}

We recall that the test-functions satisfy $\curl\Phi^\eta = \Psi^\eta$
in $\Omega_\eta$. This implies that the second integral in \eqref
{eq:testedH-hc} coincides with the left-hand side of \eqref
{eq:testedE-hc} up to the factor $i\omega\mu$ and the small error term
\begin{equation}
  \label{eq:I0eta}
  I_{0}^\eta \coloneqq \int_{\Sigma_\eta} H^\eta\cdot \vth \Psi^\eta\,.
\end{equation}
We note already here that this term vanishes in the limit $\eta\to 0$.
The volume of $\Sigma_\eta$ is of the order
$\eta^2\delta_\eta^2 / \eta = \eta\delta_\eta^2$. Inserting the bound
$|\Psi^\eta| \le C/\delta_\eta$, see the comments after \eqref
{eq:cell-functions-PP-extended}, we find
\begin{equation}
  \label{eq:I0}
  \int_{\Sigma_\eta} |\Psi^\eta|^2 \le C \int_{\Sigma_\eta}
  \left| \frac1{\delta_\eta}\right|^2
  \le C \eta\delta_\eta^2/\delta_\eta^2 = C\eta\to 0\,.
\end{equation}
Since $H^\eta$ is bounded in $L^2(\Omega, \C^3)$ by assumption, we
conclude that $I_{0}^\eta$ of \eqref {eq:I0eta} vanishes in the limit
$\eta\to 0$.

We combine the above equations \eqref {eq:testedH-hc}--\eqref
{eq:testedE-hc} into
\begin{equation}
  \label{eq:combined-EH-h}
  I_{0}^\eta 
  - \int_\Omega i \omega  \eps_\eta E^\eta\cdot \vth \Phi^\eta
  + \int_\Omega \He \cdot \Phi^\eta \times \nabla\vth
  = \frac{i}{\omega\mu} \int_\Omega \Ee \cdot \Psi^\eta \times \nabla\vth
  - \frac{i}{\omega\mu} \int_{\Sigma_\eta} \vth\, \Ee \cdot \curl \Psi^\eta\,.
\end{equation}
Denoting the integrals by $I_{1}^\eta$ to $I_{4}^\eta$, we write this
relation as
$I_{0}^\eta + I_{1}^\eta + I_{2}^\eta = I_{3}^\eta + I_{4}^\eta$.

\begin{lemma}[The limits of $I_{1}^\eta$ and $I_{2}^\eta$]
  \label{lem:easy-1h}
  There holds, as $\eta\to 0$,
  \begin{equation}
    \label{eq:easy-1h}
    I_{1}^\eta = -\int_\Omega i \omega \eps_\eta E^\eta\cdot \vth \Phi^\eta
    \to -\int_\Omega i \omega\eps_0 E\cdot \vth \,(\alpha + |x_3|)\, e_1\,,
  \end{equation}
  and
  \begin{equation}
    I_{2}^\eta = \int_\Omega \He \cdot \Phi^\eta \times \nabla\vth
    \to \int_\Omega H \cdot (\alpha + |x_3|)\, e_1 \times \nabla\vth\,.
    \label{eq:easy-2h}
  \end{equation}
\end{lemma}

\begin{proof}
  We recall that $\Phi^\eta$ is vanishing in $\Sigma_\eta$; we can
  therefore replace the factor $\eps_\eta$ by $\eps_0$. For the
  convergence, it is sufficient to exploit the strong
  $L^2$-convergence $\Phi^\eta \to (\alpha + |x_3|)\, e_1$.
\end{proof}

\begin{lemma}[The limit of $I_{3}^\eta$]
  \label{lem:easy-3high}
  There holds, as $\eta\to 0$,
  \begin{equation}
    \label{eq:easy-3high}
    (-i\omega\mu)\,  I_{3}^\eta 
    = \int_\Omega \Ee \cdot \Psi^\eta \times \nabla\vth
    \to \int_\Omega E \cdot \sign(x_3) e_2 \times \nabla\vth
    - \frac{\alpha}{\eps_0} \int_\Gamma j\, \del_1^2 \vth\,.
  \end{equation}
\end{lemma}

\begin{proof}
  In contrast to the proof for perfect conductors, where we extended
  $\Psi^\eta$ trivially to the wires, the left-hand side of \eqref
  {eq:easy-3high} contains a contribution from the wire integrals. We
  estimate this contribution with
  \begin{align*}
    \left| \int_{\Sigma_\eta} \Ee \cdot \Psi^\eta \times \nabla\vth \right|^2
    \le C \| \Ee \|^2_{L^2(\Omega)}\, \int_{\Sigma_\eta} |\Psi^\eta|^2
    \to 0\,,
  \end{align*}
  where the convergence was already derived in \eqref {eq:I0}.

  The rest of the proof is verbatim as in the case of perfect
  conductors, see Lemma~\ref {lem:easy-3p}.
\end{proof}

\begin{lemma}[The limit of $I_{4}^\eta$]
  \label{lem:easy-4h}
  In the limit $\eta\to 0$, there holds
  \begin{equation}
    \label{eq:easy-4h}
    (i\omega\mu)\,  I_{4}^\eta = \int_{\Sigma_\eta} \vth\, \Ee \cdot \curl \Psi^\eta
    \to \frac{2i}{\beta\kappa} \int_\Gamma  \vth\, j\,.
  \end{equation}
\end{lemma}

\begin{proof}
  We recall that $\Psi^\eta$ is extended into the wires and does not
  vanish in $\Sigma_\eta$. In the calculation below, we use, in this
  order: (i) $\eps_\eta = \eps_0 + \frac{i\kappa}{\theta_\eta}$ inside
  the wires, compare \eqref{eq:eps-eta}. (ii)
  $\curl \Psi^\eta = - 2/|\cC_\eta| \eta^{-1}\, e_1$ in $\Sigma_\eta$,
  see \eqref {eq:Curl-Psi}.  (iii) The choice
  $\theta_\eta = \eta\delta_\eta^2$ from \eqref {eq:delta-theta}
  together with the fact that the cross-sectional volume parameter is
  $\beta = \lim_{\eta\to 0} (|\cC_\eta|/\delta_\eta^2) =
  |\cC_{2D}|$. (iv) The flux quantity $J = j$ satisfies \eqref
  {eq:J-def}, see Lemma \ref {lem:Jj}.
  \begin{align*}
    &\int_{\Sigma_\eta} \vth\, \Ee \cdot \curl \Psi^\eta
      =  \int_{\Sigma_\eta} \vth\, \eps_\eta \Ee \cdot \frac{\theta_\eta}{i\kappa}
      \curl \Psi^\eta
    + o(1)\\
    &\qquad = 
      - \int_{\Sigma_\eta} \vth\, \eps_\eta \Ee_1 \cdot
      \frac{\theta_\eta}{i\kappa} \frac{2}{\eta |\cC_\eta|}
    + o(1)
    \to
    \frac{2i}{\beta \kappa} \lim_{\eta\to 0} \int_{\Sigma_\eta} \vth\, \eps_\eta \Ee_1
      = \frac{2i}{\beta\kappa} \int_\Gamma  \vth\, j\,.
  \end{align*}
  This was the claim.
\end{proof}

The above lemmas allow us to take the limit in the relation
$I_{0}^\eta + I_{1}^\eta + I_{2}^\eta = I_{3}^\eta + I_{4}^\eta$ of
\eqref {eq:combined-EH-h}. We find
\begin{equation}
  \label{eq:res-high-cond}
  \begin{split}
  &-\int_\Omega i \omega\eps_0 E\cdot \vth \,(\alpha + |x_3|)\, e_1
  + \int_\Omega H \cdot (\alpha + |x_3|)\, e_1 \times \nabla\vth\\
  &\qquad = \frac{i}{\omega \mu} \int_\Omega E \cdot \sign(x_3) e_2 \times \nabla\vth
  - \frac{i\alpha}{\omega \mu\eps_0} \int_\Gamma j\, \del_1^2 \vth
  +  \frac{1}{i \omega \mu}  \frac{2i}{\beta\kappa} \int_\Gamma  \vth\, j\,.
\end{split}
\end{equation}

The bulk terms have well-known contributions, we recognize
$\curl H = -i \omega \eps_0 E$ tested with
$\vth\, (\alpha + |x_3|)\, e_1$ and $\curl E = i \omega \mu H$ tested
with $(i\omega\mu)^{-1} \vth\, \sign(x_3) e_2$.  The corresponding
integrations by parts introduce boundary terms. We can extract
\begin{equation}
  \label{eq:res-high-cond-Gamma}
  \int_\Gamma \llbracket H_2 \rrbracket_\Gamma\,\alpha\vth
  = -\frac{i}{\omega \mu} \int_\Gamma E_1\, 2  \vth
  - \frac{i\alpha}{\omega \mu\eps_0}  \int_\Gamma j\, \del_1^2 \vth
  + \frac{2}{\omega \mu\beta\kappa} \int_\Gamma  \vth\, j\,.
\end{equation}
Inserting $\llbracket H_2 \rrbracket_\Gamma = i\omega j$ in \eqref
{eq:res-high-cond-Gamma}, multiplying with $i\omega\mu\eps_0/\alpha$
and using that $\vth$ was arbitrary, we arrive at
\begin{equation}
  \label{eq:bdry-cond-hc}
  -\omega^2 \mu\eps_0\, j
  - \del_1^2\, j - \frac{2i \eps_0}{\alpha\beta\kappa}\, j
  = \frac{2\eps_0}{\alpha} E_1|_\Gamma\,.
\end{equation}
The result coincides with our claim in \eqref {eq:int-law-hc}.

\appendix

\section{Properties of the cell functions}
\label {app.UV}

We study the two-dimensional cylinder
$\cZ \coloneqq (-1/2,1/2)\times \R$ with the coordinates
$(y_2,y_3)$. In order to simplify the notation, we assume that the
obstacle is centered in the origin. We use
$\cC_{2D}^\eta \coloneqq \delta_\eta\, \cC_{2D} \subset \cZ$ and
recall that the shape of the inclusion is given by a simply connected
closed set $\cC_{2D}\subset \bar B_1(0)\subset [-1,1]^2$ with a
$C^{2,\gamma}$-boundary curve.  The function
$u_\eta \colon \cZ\to \R$, $u_\eta = u_\eta(y_2,y_3)$, is defined as
the unique solution of the problem
\begin{align*}
  \Delta u_\eta &= 0 \text{ on } \cZ\setminus \cC_{2D}^\eta\,,\\
  u_\eta(y) &= 0\text{ on } \del \cC_{2D}^\eta\,,\\
  \nabla u_\eta(y) - \eta e_3\, \sign(y_3)
                &  \text{ is of class  } L^2(\cZ \setminus \cC_{2D}^\eta)\,,\\
  u_\eta(y) &\text{ periodic on } \{y_2=-1/2\}\text { and } \{y_2=1/2\}\,.
\end{align*}

The constructions that we use have similarities to those of \cite
{Schweizer-Neumann-sieve-2020} and we refer to that reference for some
points that are only sketched here. A maximum principle shows the
non-negativity of solutions, $u_\eta\ge 0$.

We recall that the inclusion has a very small diameter, we use
$\delta_\eta = e^{-\alpha\pi/\eta}$.  

\begin{lemma}[Existence, uniqueness, and a representation formula]
  For fixed $\eta>0$, the above problem has a unique solution
  $u = u_\eta$. For $y_3 > 1$, the solution can be written with
  Fourier coefficients $c_j \in \C$ as
  \begin{equation}
    \label{eq:Fourier-series}
    u(y_2,y_3) = c_0 + \eta (y_3 - 1)
    + \sum_{j\in\Z \setminus \{0\}} c_j e^{2\pi ij y_2} e^{- 2\pi |j| (y_3-1)} \,.
  \end{equation}
  An analogous representation holds for $y_3 < -1$ with Fourier
  coefficients $c_j^-$.
\end{lemma}

\begin{proof} The proof uses only standard concepts.

  \smallskip {\em Step 1: Existence.}  We choose a smooth function
  $H_\eta:\R\to \R$ satisfying $H_\eta|_{[-1,1]} \equiv 0$ and
  $H_\eta(y_3) = \eta\, |y_3|$ for every $|y_3| > 2$. The idea is to
  construct the solution $u$ in the form
  $u(y_2,y_3) = H_\eta(y_3) + w(y_2,y_3)$. We seek a function $w$
  in the space $\dot H^1_\#(\cZ\setminus \cC_{2D}^\eta)$ of
  $L^2_\loc(\cZ\setminus \cC_{2D}^\eta)$ functions with gradient in
  $L^2(\cZ\setminus \cC_{2D}^\eta)$. Additionally, we include the
  periodicity condition and the Dirichlet condition in the space.

  We would like to solve the problem
  \begin{equation}
    \label{eq:Lax-M1}
    \int_{\cZ\setminus \cC_{2D}^\eta} \nabla (H_\eta + w)\cdot\nabla\fhi = 0
  \end{equation}
  for all $\fhi$, but there is not a good class of test-functions that
  could be used. Note that $\nabla H_\eta$ is not of class
  $L^2(\cZ\setminus \cC_{2D}^\eta)$.  We therefore reformulate the
  task and consider the equation
  \begin{equation}
    \label{eq:Lax-M2}
    \int_{\cZ\setminus \cC_{2D}^\eta} \nabla w\cdot\nabla\fhi
    =  \int_{\cZ\setminus \cC_{2D}^\eta} (\Delta H_\eta)\ \fhi
    \qquad
    \forall \fhi\in \dot H^1_\#(\cZ\setminus \cC_{2D}^\eta)\,.
  \end{equation}
  In this form, the problem can be solved with the lemma of
  Lax--Milgram; we observe that $\Delta H_\eta$ has a bounded support.
  The test-function $\fhi = w$ yields an $L^2$-bound for the gradient,
  namely
  \begin{equation}
    \label {eq:grad-u-bound}
    \int_{\cZ\setminus \cC_{2D}^\eta} |\nabla w|^2
    \le C_P(\eta)\, \| \Delta H_\eta \|_{L^2(\cZ\setminus \cC_{2D}^\eta)}^2\,,
  \end{equation}
  where $C_P(\eta)$ is a constant from the Poincar\'e estimate; it
  deteriorates with $\eta\to 0$ (in the sense that
  $C_P(\eta)\to \infty$ as $\eta\to 0$) since the Dirichlet boundary
  gets smaller.  We remark that we can choose $H_\eta$ with
  $\| \Delta H_\eta \|_{L^2(\cZ\setminus \cC_{2D}^\eta)} \le C \eta$.
  
  \smallskip {\em Step 2: Uniqueness.} Let $u^{(1)}$ and $u^{(2)}$ be
  two solutions.  The difference $\psi = u^{(1)} - u^{(2)}$ is also a
  harmonic function, vanishing on the inner boundary, but it has
  another growth property: $\psi$ satisfies
  $\nabla \psi \in L^2(\cZ\setminus \cC_{2D}^\eta)$.  We can use
  $\psi$ as a test-function for $\Delta \psi = 0$ and obtain that
  $\nabla \psi$ vanishes identically (integration by parts is indeed
  possible, this must be verified using that averages of $\psi$ over
  lines $\{ y_3 = \xi \}$ are independent of $\xi>1$). As a
  homogeneous Dirichlet condition is imposed on a part of the
  boundary, we conclude $\psi = 0$ and hence uniqueness.

  \smallskip {\em Step 3: Representation formula.} We consider a
  solution $u$ and recall that, in this proof, the parameter $\eta$ is
  fixed. As a harmonic function that satisfies a Dirichlet condition
  along $\del \cC_{2D}^\eta$, the solution is, on the set
  $[(-1/2,1/2)\times (-2,2)]\setminus \cC_{2D}^\eta$, a function of class
  $H^1$, along the lines $\{ y_3 = 1\}$ and $\{ y_3 = -1\}$, it is of
  class $L^2$. This allows us to expand $u|_{\{ y_3 = 1\}}$ in a
  Fourier-series with coefficients $c_j$, the complex number $c_0$ is
  the average of $u$ on the line $\{ y_3 = 1\}$.

  When we consider the corresponding extension $\tilde u$ as given by
  the expansion \eqref {eq:Fourier-series}, we observe: (1) $\tilde u$
  coincides with $u$ on the line. (2) $\tilde u$ is harmonic and
  periodic in $y_2$. (3) There holds
  $\nabla \tilde u - \eta\sign(y_3) e_3 \in L^2(\cZ\setminus
  \cC_{2D}^\eta)$.  The function $u$ has the same properties (1)--(3)
  on the set $\{ y_3 > 1 \}$. As these properties determine the
  solution uniquely, we find $u = \tilde u$ on $\{ y_3 >1 \}$ and have
  thus verified the representation formula. It shows, in particular,
  that $u(y_2,y_3) - \eta\, |y_3|$ is bounded.
\end{proof}

We note that the above lemma provides also one other claim about the
shape function $U_1^\eta$, namely
$U_1^\eta(y) - \eta y_3 \to c_0 - \eta$ for $y_3 \to \infty$; we
recall that $c_0$ also depends on $\eta$.

\medskip In the next lemma appears the parameter $\alpha$ of
averages. We emphasize that it is the parameter that was used to
define the small obstacle sizes, $\delta_\eta = e^{-\alpha\pi/\eta}$
in \eqref {eq:delta-theta}.

\begin{lemma}[Estimates for $u_\eta$]
  There exists a constant $C$, independent of $\eta$, such that
  \begin{align}
    \label{eq:uni-est-ueta-1}
    \| u_\eta |_{\{ y_3 = \pm 1 \}} \|_{L^2(\{ y_3 = \pm 1 \})}
    &\le C\,.
  \end{align}
  The averages of $u_\eta$ on $\{ y_3 = \pm 1 \}$ satisfy
  \begin{equation}
    \label{eq:alpha-eta-conv}
    c_{0,\eta}^+ \to \alpha\quad \text{ and }\quad
    c_{0,\eta}^- \to \alpha\,.
  \end{equation}

  We extend $u_\eta$ periodically in $y_2$-direction and extend it by
  zero into the obstacles. Using this extension, we define
  $\phi_\eta \colon (0,1)\times (-1,1) \to \R$ by rescaling, setting
  $\phi_\eta(x) \coloneqq u_\eta(x/\eta)$. Then $\phi_\eta$ satisfies
  \begin{equation}
    \phi_\eta \to \phi_0\text{ in } L^2((0,1)\times (-1,1))\,,
    \quad \text{ where }\quad
    \phi_0(x) \coloneqq \alpha + |x_3|\,.
    \label{eq:uni-est-ueta-3}
  \end{equation}
\end{lemma}

\begin{proof}
  \smallskip {\em Step 1: Averages on $B_{1/2}(0)$.} This step is the
  core of the proof. It is based on the Theorem of Gauß. We exploit
  that, on circles $\del B_{r}(0)$ with radius $r\le 1/2$ and outer
  normal vector $\nu$, the ``fluxes'' are independent of $r$. More
  precisely: As $u_\eta$ is harmonic between two such circles, we
  obtain that the flux quantities
  $\int_{\del B_{r}(0)} \nu\cdot\nabla u_\eta$ are independent of $r$.
  With the same argument, this quantity coincides with
  $\int_{\{ y_3 = +h\}} e_3\cdot\nabla u_\eta + \int_{\{ y_3 = -h\}}
  (-e_3)\cdot\nabla u_\eta = 2\eta$.
  
  This allows us to calculate the derivatives of averages:
  \begin{align*}
    &\del_r \left[ \frac1{2\pi r} \int_{\del B_{r}(0)} u_\eta(y)\, dS(y) \right]
      = \del_r \left[ \frac1{2\pi} \int_{\del B_{1}(0)} u_\eta(r y)\, dS(y) \right]\\
    &\qquad= \frac1{2\pi} \int_{\del B_{1}(0)} \nu\cdot\nabla u_\eta(r y)\, dS(y)
      = \frac1{2\pi r} \int_{\del B_{r}(0)} \nu\cdot\nabla u_\eta(y)\, dS(y)
      = \frac{\eta}{\pi r}\,.
  \end{align*}

  In order to exploit this information, we must consider the average
  over a small ball. We recall that the closed ball of radius
  $\delta_\eta$ contains the inclusion $\cC_{2D}^\eta$. We denote the
  average of $u_\eta$ over this smallest circle by $A_\eta$,
  \begin{equation}
    A_\eta := \frac1{2\pi \delta_\eta}
    \int_{\del B_{\delta_\eta}(0)} u_\eta(y)\, dS(y)\,.
  \end{equation}
  Using this quantity, we can proceed with the above argument.  We
  integrate the above derivatives of averages from the radius
  $\delta_\eta$ to the radius $1/2$ and find
  \begin{align}
    \nonumber
    \alpha_\eta
    &\coloneqq \frac1{\pi} \int_{\del B_{1/2}(0)} u_\eta
      = A_\eta  + \int_{\delta_\eta}^{1/2} \frac{\eta}{\pi r}\, dr
      = A_\eta  + \left. \frac{\eta}{\pi}\ \log(r)\right|_{\delta_\eta}^{1/2}\\
    &= A_\eta + \frac{\eta}{\pi}\ \log(1/2)
      - \frac{\eta}{\pi}\ \log(\delta_\eta)
      = A_\eta + \frac{\eta}{\pi}\, \log(1/\delta_\eta)
      + O(\eta)
      \nonumber\\
    &= A_\eta + \frac{\eta}{\pi}\, \log(e^{\alpha\pi/\eta}) 
      + O(\eta)
      = A_\eta + \alpha + O(\eta) \,.
      \label{eq:average-on-B12}
  \end{align}

  \smallskip {\em Step 2: Vicinity of the inclusion and a gradient
    estimate.}  Estimate \eqref {eq:average-on-B12} provides, in
  particular, that the average values of $u_\eta$ are bounded by
  $A_\eta + C$ on $\del B_{1/2}(0)$ for some $\eta$-independent
  constant $C$. As the function $u_\eta$ is harmonic, non-negative,
  and has bounded gradients for $|y_3|\to \infty$, on bounded sets,
  $u_\eta$ is pointwise bounded by $C (A_\eta + 1)$. This can be
  exploited in the Lax-Milgram description \eqref {eq:Lax-M2}:
  Inserting $\fhi := u_\eta - H_\eta$, we find that the right-hand
  side is bounded by $C \eta (A_\eta + 1)$ on bounded subsets and we
  conclude, for $\cZ_1 \coloneqq (-1/2,1/2)\times (-1,1)$,
  \begin{equation}
    \label{eq:sqrt-eta-grad-est}
    \int_{\cZ_1 \setminus \cC_{2D}^\eta} |\nabla u_\eta|^2
    \le C \eta (A_\eta + 1)\,.
  \end{equation}
  We use a Poincar\'e estimate on the smallest ball, using the
  homogeneous Dirichlet condition on $\del\cC_{2D}^\eta$, to find
  \begin{equation*}
    A_\eta \le C\, 
    \left( \int_{B_{\delta_\eta}(0) \setminus \cC_{2D}^\eta}
      |\nabla u_\eta|^2\right)^{1/2}
    \le C \sqrt{\eta} \sqrt{A_\eta + 1}\,.
  \end{equation*}
  This provides $A_\eta \le C \sqrt{\eta}$.
  
  \smallskip {\em Step 3: Uniform bounds and convergence of averages.}
  We have already observed in Step 2 that, on bounded subsets,
  $u_\eta$ is bounded by $C (A_\eta + 1)$. As $A_\eta$ is bounded and
  $u_\eta$ is non-negative, this provides also the $L^2$-estimate on
  the horizontal interfaces as stated in \eqref {eq:uni-est-ueta-1}.

  Using \eqref {eq:sqrt-eta-grad-est} and a Poincar\'e estimate
  outside the ball $B_{1/2}(0)$, we also see that averages on
  different parts of the boundary can only differ in the order
  $\sqrt{\eta}$,
  \begin{equation}
    \label{eq:averages-pm}
    | c_{0,\eta}^+ - \alpha_\eta | =
    \left| \int_{\{ y_3 = +1\}} u_\eta - \alpha_\eta \right|
    \le C \sqrt{\eta}\,,
  \end{equation}
  and likewise for $-$ instead of $+$.  Together with \eqref
  {eq:average-on-B12}, we obtain \eqref {eq:alpha-eta-conv}.

  \smallskip {\em Step 4: Strong convergence of rescaled functions.}
  We can calculate the convergence directly with the help of the
  representation formula. We use the set
  $R = \{ (x_2,x_3)\in (0,1)^2\,|\,x_3 > \eta\}$ and insert \eqref
  {eq:Fourier-series} in the second equality:
  \begin{align*}
    &\int_R | \phi_\eta - \phi_0 |^2
      = \int_R | u_\eta(x_2/\eta,x_3/\eta) - \alpha - |x_3| |^2\, dx_2\, dx_3\\
    &\qquad = \int_R \left| c_{0,\eta}^+ + x_3 - \eta - \alpha - |x_3|
      + \sum_{j\in\Z \setminus \{0\}}
      c_j e^{2\pi ij x_2/\eta} e^{- 2\pi |j| (x_3/\eta - 1)} \right|^2\, dx_2\, dx_3\\
    &\qquad \le C \Big( \eta + |c_{0,\eta}^+ - \alpha|^2 
      + \int_R \sum_{j\in\Z \setminus \{0\}}
      \left| c_j e^{2\pi ij x_2/\eta} e^{- 2\pi |j| (x_3/\eta - 1)} \right|^2\,
      dx_2\, dx_3 \Big)\\
    &\qquad \le C \Big( \eta + |c_{0,\eta}^+ - \alpha|^2 
      + \sum_{j\in\Z \setminus \{0\}} |c_j|^2  \int_\eta^1 \left| e^{- 2\pi (x_3/\eta - 1)}
      \right|^2\, dx_3 \Big)
      \le\, C\, \eta\,.
  \end{align*}
  The calculation is done similarly in the lower part and the integral
  over the strip $-\eta < x_3 < \eta$ is also of order $\eta$ by the
  uniform boundedness of $u_\eta$. This provides \eqref
  {eq:uni-est-ueta-3}.
\end{proof}

\medskip It remains to discuss the function
$V^\eta = \eta^{-1} \curl(U^\eta)$ that was defined in \eqref
{eq:V-eta}. In our two-dimensional setting, we have to study the
gradient of $u_\eta$.  The stabilization property (1) of the
subsequent lemma provides \eqref {eq:V-stabilization}. The bound of
(2) is used only for the extension of $V^\eta$ into the wire in \eqref
{eq:V-eta-cC}. More precisely: The function $g$ has typical values of
order $\eta/\delta_\eta$, hence the gradient of $U_1^\eta$ has also
typical values of order $\eta/\delta_\eta$ and $V^\eta$ has typical
values of order $1/\delta_\eta$.

\begin{lemma}[Properties of $v_\eta$] The function
  $v_\eta \coloneqq \eta^{-1} \nabla u_\eta$ has the following
  properties:\smallskip

  (1) Stabilization: $v_\eta(y) \to \pm e_3$ as $|y_3|\to \infty$.

  (2) Uniform bound: There holds $|v_\eta| \le
  C/\delta_\eta$. Furthermore, rescaling the domain by a factor
  $\delta_\eta$, the function $v_\eta$ and its derivative are bounded
  on $\del\cC_{2D}$ also by $C/\delta_\eta$.
\end{lemma}

Before starting the proof, we observe the following: In the vicinity
of the inclusion, we know that the flux quantity
$\int_{\del B_{r}(0)} \nu\cdot\nabla u_\eta$ is $2\eta$, independent
of $r$, hence
\begin{equation}
  \label{eq:flux-v}
  \int_{\del B_{r}(0)} \nu\cdot v_\eta = 2\,.
\end{equation}
This shows that $v_\eta$ must have at least values of the order
$|v_\eta| \sim 1/\delta_\eta$ near the obstacle.

\begin{proof}
  The stabilization property (1) follows immediately from the
  representation formula \eqref {eq:Fourier-series} of $u_\eta$.

  \smallskip Regarding (2) we note that $\nabla u_\eta$ is at most of
  the order $\eta$ for $|y_3|\ge 2$, this provides the bound
  $|v_\eta| \le C$ in this part of the domain.  The uniform bound on
  $\cZ_2\setminus B_{\delta_\eta}$ for
  $\cZ_2 \coloneqq (-1/2,1/2)\times (-2,2)$ follows by a maximum
  principle from the above when we show the uniform bound in the
  vicinity of $\del B_{\delta_\eta}$.

  It remains to show the bound in the vicinity of
  $\del B_{\delta_\eta}$. We consider a rescaled situation and the
  domain $B_2(0)\setminus\cC_{2D}$.  Let $g$ be the harmonic function
  on this domain with $g = 0$ on $\del\cC_{2D}$ and $g = 1$ on
  $\del B_2(0)$. By the strong maximum principle and the Hopf lemma,
  there holds $c_0 \coloneqq \min_{\del B_2(0)} \nu\cdot\nabla g >
  0$. We can use the rescaled functions
  $g_\eta(y) \coloneqq g(y/\delta_\eta)$ as test-functions in the
  equation $\Delta u_\eta = 0$ on
  $B_{2\delta_\eta}(0)\setminus\cC_{2D}^\eta$ and find
  \begin{equation}
    \label{eq:Green-ball}
    \frac{c_0}{\delta_\eta} \int_{\del B_{2\delta_\eta}(0)} u_\eta
    \le \int_{\del B_{2\delta_\eta}(0)} u_\eta\, \nu\cdot\nabla g_\eta
    = \int_{\del B_{2\delta_\eta}(0)} g_\eta\, \nu\cdot\nabla u_\eta
    = \int_{\del B_{2\delta_\eta}(0)} \nu\cdot\nabla u_\eta
    = 2\eta\,.
  \end{equation}

  We set $h_\eta(z) \coloneqq \eta^{-1} u_\eta(\delta_\eta z)$. The
  function $h_\eta$ is non-negative and harmonic on
  $B_3(0)\setminus\cC_{2D}$, it vanishes on $\del\cC_{2D}$, and
  \eqref{eq:Green-ball} provides
  $\int_{\del B_2(0)} h_\eta \le 2/c_0$.  The Harnack inequality on
  the annulus $\{3/2 < |z| < 3\}$ yields an $\eta$-independent bound
  for $h_\eta$ on $\del B_2(0)$, and the maximum principle yields the
  same bound on $B_2(0)\setminus\cC_{2D}$. Since $\del\cC_{2D}$ is of
  class $C^{2,\gamma}$, elliptic regularity up to the boundary gives
  $|\nabla h_\eta| \le C$ on $B_{3/2}(0)\setminus\cC_{2D}$. Because of
  $\nabla h_\eta(z) = \delta_\eta\, v_\eta(\delta_\eta z)$, this is
  the bound $|v_\eta| \le C/\delta_\eta$ on
  $B_{3\delta_\eta/2}(0)\setminus\cC_{2D}^\eta$, in particular in the
  vicinity of $\del B_{\delta_\eta}$. As the boundary is of class
  $C^{2,\gamma}$, Schauder estimates also provide bounds for second
  derivatives of $h_\eta$ on $\del\cC_{2D}$. This provides (2).
\end{proof}

\section{An $L^1(\Omega)$-estimate and the limiting flux measure}
\label {app.limit-J}

In this short appendix we draw some conclusions from the
$L^2$-estimate for solutions in \eqref {eq:uniform-L2-bound} in the
highly conducting case. They are quite standard and parallel to, e.g.,
\cite {Bouchitte-Schweizer-Max-2010}.

A first step is to use the solutions as test-functions in the
Maxwell-system in order to obtain improved estimates. More precisely,
we multiply \eqref {eq:Strong:Maxwell:eta:1} with $\bar H^\eta$ and
the complex conjugate of \eqref {eq:Strong:Maxwell:eta:2} with
$E^\eta$.  Integrating produces two identical term on the left-hand
side and we find uniform $L^1(\Omega)$-bounds for
$|\eps_\eta|\, |E^\eta|^2$ and $|\mu|\, |H^\eta|^2$.  This estimate
allows to bound integrals over wires for our choice of $\eps_\eta$ and
the geometry:
\begin{align*}
  \int_{\Sigma_\eta} |\eps_\eta|\ |E^\eta|
  &\le \left( \int_{\Sigma_\eta} |\eps_\eta|\ |E^\eta|^2 \right)^{1/2}
  \left( \int_{\Sigma_\eta} |\eps_\eta| \right)^{1/2}\\
  &\le C (|\Sigma_\eta|\ \theta_\eta^{-1})^{1/2}
  \le C ((\eta\delta_\eta^2)/(\eta\delta_\eta^2))^{1/2} = C\,.
\end{align*}
We obtain that the flux quantity $F^\eta \coloneqq \eps_\eta E^\eta$
is bounded in $L^1(\Omega)$. Hence, up to a subsequence, it possesses
a limiting measure. By the $L^2$-boundedness of $E^\eta$, the singular
part $\vec J$ of the limiting measure, defined in \eqref {eq:J-def},
must be concentrated on $\Gamma$.

It is easy to show $\vec J = J\, e_1$, that $\vec J$ points in
direction $e_1$. We present the argument for $\vec J\cdot e_2 = 0$,
the argument for $e_3$ follows with the analogous calculation. We
choose a smooth periodic test-function $q : [0,1]\to \R$ with
$q'\equiv 1$ in a central part (at the $y_2$-position of the
wire). With an arbitrary smooth and compactly supported function
$\vth:\Omega\to \R$, we use the test-function
$\fhi(x) := \eta \nabla [q(x_2/\eta)\, \vth(x)]$ in \eqref
{eq:Strong:Maxwell:eta:2}. The left-hand side vanishes as the
test-function is a gradient and the right-hand side has the limit
$-i\omega \int_\Gamma \vec J\cdot e_2\, \vth$.

\printbibliography

\end{document}